%% file: main.tex
\documentclass[11pt]{amsart}
\usepackage[total={6.5in,8.75in}, top=1.2in, left=0.9in, includefoot]{geometry}

\usepackage{cleveref}
\usepackage{algorithm}
\usepackage{url}
\usepackage[symbol]{footmisc}

\newcommand{\lsp}{\vspace{3mm}}

\theoremstyle{definition}

\numberwithin{definition}{section}

\usepackage[backend=biber]{biblatex}
\input{shared}

\input{notations}

\begin{document}

\begin{center}

\textsc{A Simplified Fast Multipole Method\\ Based on Strong Recursive Skeletonization}

\lsp

Anna Yesypenko\footnote{\label{note:oden}Oden Institute, University of Texas at Austin (annayesy@utexas.edu, pgm@oden.utexas.edu).}, 
Chao Chen\footnote{Department of Mathematics, North Carolina State University  (cchen49@ncsu.edu).}, 
and Per-Gunnar Martinsson\textsuperscript{\ref{note:oden}}\\

\lsp
\end{center}

\begin{center}
\begin{minipage}{130mm}
\small
\textbf{Abstract:}
\input{abstract}
\end{minipage}
\end{center}

\input{skelfmm_paper}

\printbibliography

\end{document}

%% file: shared.tex
\usepackage{amssymb}
\usepackage[shortlabels]{enumitem}
\usepackage{algorithmicx}
\usepackage[noend]{algpseudocode}
\usepackage{bm}
\usepackage{booktabs}
\usepackage{capt-of}
\usepackage{color, colortbl}
\usepackage{mathtools}
\usepackage{multirow}
\usepackage{siunitx}
\usepackage{subcaption}

\usepackage{xcolor}

\usepackage{booktabs}
\usepackage{multirow}
\usepackage{siunitx}
\usepackage{tikz}

\usepackage{lipsum,pgfplots}
\usepackage{amsfonts}
\usepackage{graphicx}
\usepackage{epstopdf}
\usepackage{amsopn}
\usepackage{changes}

\ifpdf
  \DeclareGraphicsExtensions{.eps,.pdf,.png,.jpg}
\else
  \DeclareGraphicsExtensions{.eps}
\fi

\graphicspath{{figs/}}

%% file: notations.tex
\newcommand{\vct}[1]{\bm{\mathsf{#1}}}
\newcommand{\mtx}[1]{\bm{\mathsf{#1}}}
\newcolumntype{L}{>{$}c<{$}} 
\newcolumntype{C}{>{$}c<{$}} 

\newcommand{\A}{\mtx{A}}
\newcommand{\x}{\vct{x}}

\newcommand{\addstretch}[1]{\addtolength{#1}{\fill}}

\newcommand{\N}{ {\mathcal{N}} }

\newcommand{\bigO}{ {\mathcal{O}} }
\newcommand{\B}{ {\mathcal{B}} }

\newcommand{\X}{ {\mathcal{X}} }
\newcommand{\F}{ {\mathcal{F}} }
\renewcommand{\L}{ {\mathcal{L}} }
\newcommand{\T}{ {\mathcal{T}} }
\newcommand{\eps}{ {\varepsilon} }

\newcommand{\R}{ {\mathcal{R}} }
\newcommand{\C}{ {\mathcal{C}} }

\renewcommand{\P}{ {\mathcal{P}} }
\renewcommand{\S}{ {\mathcal{S}} }

\newcolumntype{H}{>{\setbox0=\hbox\bgroup}c<{\egroup}@{}}

\newcommand{\pluseq}{\mathrel{+}=}
\newcommand{\minuseq}{\mathrel{-}=}

%% file: abstract.tex
This work introduces a kernel-independent, multilevel, adaptive algorithm for efficiently evaluating a discrete convolution kernel with a given source distribution.
The method is based on linear algebraic tools such as low rank approximation and ``skeleton representations'' to approximate far-field interactions. 
While this work is related to previous linear algebraic formulations of the fast multipole method, the proposed algorithm is distinguished by relying on simpler data structures.

\vspace{1.5mm}

The proposed algorithm eliminates the need for explicit interaction lists by restructuring computations to operate exclusively on the near-neighbor list at each level of the tree, thereby simplifying both implementation and data structures. This work also introduces novel translation operators that significantly simplify the handling of adaptive point distributions.
As a kernel-independent approach, it only requires evaluation of the kernel function, making it easily adaptable to a variety of kernels.  By using operations on the neighbor list (of size at most 27 in 3D) rather than the interaction list (of size up to 189 in 3D), the algorithm is particularly well-suited for parallel implementation on modern hardware.

\vspace{1.5mm}

Numerical experiments on uniform and non-uniform point distributions in 2D and 3D demonstrate the effectiveness of the proposed parallel algorithm for Laplace and (low-frequency) Helmholtz kernels. The algorithm constructs a tailored skeleton representation for the given geometry during a precomputation stage. After precomputation, the fast summation achieves high efficiency on the GPU using batched linear algebra operations.

%% file: skelfmm_paper.tex
\section{Introduction}
We present an algorithm for evaluating a sum of the form
\begin{equation}
u_i  = \sum_{j=1,i\neq j}^N G(x_i,x_j)\ q_j, \qquad i=1,\dots N
\label{e:sum_fmm}
\end{equation}
where $\X = \{x_i\}_{i=1}^N$ is a given set of points
in $\mathbb{R}^2$ or $\mathbb{R}^3$, and where $G$ is a given kernel function
associated with a standard elliptic PDE of mathematical physics.
For instance, $G$ may be the fundamental solution of Laplace's equation,
\begin{equation} \label{e:laplace}
G(x_i, x_j) = 
\left\{
\begin{array}{cc}
- \frac{1}{2\pi}\log\left(\|x_i - x_j \|\right) & \quad x_i \neq x_j \in \mathbb{R}^2  \\
\frac{1}{4\pi \|x_i - x_j \|} & \quad x_i \neq x_j \in \mathbb{R}^3.
\end{array}
\right.
\end{equation} 
The task of evaluating a sum such as (\ref{e:sum_fmm}) arises frequently in particle simulations, computational chemistry, molecular dynamics, the discretization of boundary integral equations,
and many other contexts. The summation problem (\ref{e:sum_fmm}) can be viewed as a matrix-vector product
\begin{equation}
\vct{u} = \A \, \vct{q},
\label{e:matvec}
\end{equation}
where $\A$ is a dense matrix of size $N\times N$ with off-diagonal entries $A_{ij} = G(x_i,x_j)$ and zeros on the diagonal, and where $\vct{u} = {[u_1, \dots, u_N]}$ and  $\vct{q} = {[q_1, \dots, q_N]}$.

Evaluating (\ref{e:sum_fmm}) directly requires $\mathcal O(N^2)$ operations.
Over the past decades, several analytic and algebraic algorithms have been proposed to evaluate (\ref{e:sum_fmm}) efficiently. 
For uniform particle distributions, mesh based methods leveraging the Fast Fourier Transform (FFT) achieve $\bigO(N \log N)$ complexity \cite{essmann1995smooth,toukmaji1996ewald,darden1993particle}. 
For highly adaptive point distributions, the Fast Multipole Method (FMM)
is a preferred approach, offering $\bigO(N)$ complexity by operating on a hierarchical adaptive tree that partitions the given geometric points.

The original FMM used analytic expansions to efficiently approximate far-field interactions. Specifically, multipole expansions approximate the influence of source points in the far-field, while local expansions capture the effect of far-field sources on a target region. These expansions are constructed through hierarchical tree traversals: multipole expansions are computed during the upward traversal, and local expansions are computed during the downward traversal. Once the local expansions are available in the smallest target regions (leaf boxes), the result of (\ref{e:sum_fmm}) is computed by combining the contributions from the local expansions with the contributions from near neighbors, which are evaluated directly.

The FMM was first developed for the Laplace kernel \cite{greengard1987fast,rokhlin1985rapid} and later extended to a wide range of kernels in computational physics \cite{fu1998fast,fu2000fast,yoshida2001application}, including specialized techniques for 3D applications \cite{greengard1996new,greengard1997new}. However, deriving analytic expansions for arbitrary kernels can be challenging, and implementing efficient translation operators often becomes cumbersome.

To broaden the applicability of the FMM, kernel-independent
algebraic variants \cite{ying2004kernel,malhotra2015pvfmm,martinsson2007accelerated} were introduced. 
These methods exploit the algebraic property that interactions between well-separated point sets are numerically low-rank for many non-oscillatory, translation-invariant kernel functions. Using tools such as singular value decomposition (SVD) or other numerical decompositions, these techniques compress interactions without relying on analytic formulas \cite{fong2009black,gimbutas2003generalized}.

A particularly effective approach involves selecting an equivalent set of densities, or so-called skeletons, that replicate far-field interactions
\cite{ying2004kernel, malhotra2015pvfmm,malhotra2016algorithm,martinsson2007accelerated,cheng2005compression}. 
The skeleton points are typically precomputed algebraically for a given user tolerance and the computation of the approximate far-field interactions only involves
evaluating the kernel entries,
enabling portable implementations across various architectures using standard linear algebra protocols.

While kernel-independent methods simplify the handling of arbitrary kernels, they retain the same algorithmic constructs and challenges as the original analytic FMM. Specifically, translating multipole expansions to local expansions depends on efficiently traversing auxiliary data structures such as the interaction list (V-list), W-list, and X-list \cite{ying2006kernel} to accumulate far-field interactions for each box. In 3D, the interaction list can contain up to 189 boxes, compared to only 27 boxes in the neighbor list. Navigating these larger data structures requires careful data management and load balancing, introducing significant implementation challenges, particularly for adaptive trees.

\subsection{Contributions of the Present Work}

The proposed algorithm introduces a significantly simplified kernel-independent FMM that bypasses the need for interaction lists and related auxiliary lists for adaptive trees. Instead, it relies solely on neighbor-to-neighbor interactions at each level of the tree. The key features of the method include:

\begin{itemize}
    \item Expanded framework for outgoing and incoming expansions.
    The algorithm uses four types of expansions per box (two types of outgoing expansions and two types of incoming expansions), compared to two in classical FMM (multipole and local). These additional expansions help to reorganize the algorithmic structure.

    \item Algebraic modifications to near-neighbor computations.
    Rather than directly evaluating near-neighbor contributions for the smallest target boxes, the algorithm introduces algebraically modified translations. These modifications alter the recursive nature of the multilevel approach, leading to a new set of operations needed to translate outgoing expansions to incoming expansions. These operations require traversing the near-neighbors at every level of the tree.

    \item Simplified translations for adaptive trees.
    For nonuniform trees, translations are performed on coarse neighbors (neighboring boxes at a higher level) and fine neighbors (neighboring boxes at a lower level).
\end{itemize}

Our proposed method has a minor increase in the constant prefactor
compared to the traditional FMM.
The prefactor costs of the classical FMM are a function of a desired accuracy and the size of the interaction list.
The classical FMM traverses the interaction list to translate multipole to local expansions for every box. 
There are 27 and 189 boxes in the interaction list for uniform point distributions in 2D and 3D, respectively.

Our proposed algorithm in principle traverses a slightly larger list (of sizes 36 and 216 for uniform distributions in 2D and 3D, respectively) to translate outgoing expansions to incoming expansions. However, these operations are actually accumulated at the next coarsest level of the tree, requiring instead the a traversal of the neighbor list, albeit for more skeleton points in each box.
A key advantage of the proposed approach is its reduced parallel overhead, as neighbor list operations require substantially less coordination and synchronization.

\subsection{Context and Related Work}

A common theme in recent advances in FMM-type algorithms is the improvement of practical runtime through simplified data structures, while maintaining compatibility with adaptive trees.
One such example is the dual-space multilevel kernel-splitting (DMK) framework \cite{jiang2023dual}, which operates on adaptive trees with a multilevel method inspired by Ewald’s approach \cite{ewald1921berechnung, bagge2023fast, shamshirgar2021fast}. DMK eliminates the need for the interaction list and relies solely on the near-neighbor list, though it differs fundamentally from the present work by employing analytic rather than algebraic techniques.

Direct solvers and preconditioners for systems such as (\ref{e:matvec}) --- where $\mtx A$ incorporates corrections to the diagonals to ensure invertibility --- have increasingly focused on simplified data structures to improve the practical performance of hierarchical matrix solvers \cite{hackbusch2015hierarchical,hackbusch1999sparse,borm2003introduction}. Solvers using weak admissibility compress near-neighbor interactions as numerically low-rank \cite{greengard2009fast, michielssen1996scattering, hofastdirect}, effectively eliminating the need for interaction lists, albeit with higher asymptotic complexity for 3D point distributions. Recent advances in direct solvers \cite{yesypenko2023randomized,yesypenko2023rsrs} that build on earlier works \cite{minden2017recursive, sushnikova2022fmm} construct LU or Cholesky decompositions of the FMM matrix by algebraically compressing interaction lists. This approach yields an invertible factorization that relies solely on neighbor interactions at each level of the tree.

The proposed method introduces a novel algorithmic structure that completely eliminates interaction lists from the kernel-independent FMM. By relying solely on neighbor-to-neighbor interactions, it simplifies parallel implementation on modern architectures while maintaining full compatibility with adaptive trees.

\section{Methodology}\label{s:methodology}

In this section, the summation problem and general methodology used for the FMM are discussed. As a simple example, consider the geometry shown in Figure \ref{f:1dgeom}.
The points are partitioned in a multi-level uniform tree $\mathcal T$, and the geometry is described by some simple terminology.
Boxes on the same level are called \textit{colleagues}. Colleagues that share an edge or corner (i.e.~are adjacent) are called \textit{neighbors}. Colleagues which are not adjacent are called \textit{well-separated}.
For a box $\B$, the set of neighbor boxes is called $\N(\B)$, and the \textit{far-field} $\F(\B)$ is the set of all well-separated boxes. The \textit{parent} $\P(\B)$ of a box $\B$ is the box on the next coarser level which contains $\B$. Likewise, the \textit{children} $\C(\B)$ of a box $\B$ are the set of boxes whose parent is $\B$. A box $\B$ without children is called a \textit{leaf}.
\begin{figure}[!htb]
    \centering
    \input{1d_geom}
    \caption{Points partitioned into multi-level tree $\mathcal T$.
    For a box $\B$, the set of neighbors is called $\N(\B)$, and the set of well-separated boxes is called the far field $\F(\B)$. As an example $\N(9) = \{8,10\}$ and $\F(9)=\{11,\dots,15\}$.}
    \label{f:1dgeom}
\end{figure}

\noindent The kernel matrix $\mtx A$ has the following low-rank property: if two well-separated colleague boxes $\B_i,\B_j$ each contain $m$ points, then the corresponding off-diagonal submatrix ${\mtx A_{ \B_i, \B_j}}$ satisfies
\begin{equation}
    \underset{m \times m}{\mtx A_{ \B_i, \B_j}} = \underset{m \times k}{  {\mtx E}_{ \B_{i\vphantom{j}}} }\ \underset{k \times k}{\tilde{\mtx A}_{\vphantom{\B_j}}}\ \underset{k \times m} { {\mtx F}_{\B_j}} + O(\epsilon),
    \label{e:wellsep}
\end{equation}
where ${\mtx E}_{\B_i},{\mtx F}^*_{\B_j} \in \mathbb{R}^{m \times k}$ are bases for the column and row spaces, respectively. 
In particular, the numerical rank $k \sim \bigO(\log(1/\epsilon))$ is independent of $m$, for a smooth non-oscillatory kernel  
such as the free-space Green’s function for the Laplace equation (\ref{e:laplace}) in 2D (see, e.g.,~\cite{greengard1987fast}).
For non-oscillatory kernels on 3D point distributions, the numerical rank $k_{\rm max}$ is also independent of $m$ but has slower decay with $\epsilon$, in particular, $k \sim \bigO(\log(1/\epsilon)^2)$.

The low-rank property may deteriorate for oscillatory kernels, such as the Green's function for the Helmholtz equation (parameterized by the wavenumber $\kappa >0$)
\begin{equation} \label{e:helmholtz}
G(x_i, x_j) = \begin{cases}\frac{i}{4}H_{0}^{(1)}(\kappa \|x_i - x_j \|),& x_i \neq x_j \in \mathbb{R}^2\\
\frac 1 4 \frac{e^{i \kappa \|x_i - x_j \|}}{\|x_i - x_j \|},& x_i \neq x_j \in \mathbb{R}^3,
\end{cases}
\end{equation}
where $H_0^{(1)}$ is the Hankel function of the first kind of order zero. 
In particular, for a box of diameter $D$, the numerical rank is $\bigO(\kappa D)$. 
To achieve log-linear complexity when the number of points scales with the wavenumber $\kappa$ (maintaining a constant number of points per wavelength), specialized techniques are required. These include analytic methods (e.g., \cite{rokhlin1990rapid,rokhlin1993diagonal,cheng2006wideband}) and algebraic approaches (e.g., \cite{engquist2007fast,benson2014parallel,li2015butterfly}).
When the wavenumber $\kappa$  is fixed as the number of points $m$ increases, the numerical rank remains independent of $m$. In this regime, the techniques described in this work remain applicable.

\begin{figure}[!htb]
    \centering
    \input{FMM_separable}
    \caption{The kernel matrix $\mtx A$ corresponding to geometry in Figure \ref{f:1dgeom} for the boxes on level 3
    has the algebraic structure in the figure.
    Interactions between neighbors are full rank. 
    The interaction between a box and its far field is low rank to accuracy $\epsilon$, as described in (\ref{e:wellsep}).}
    \label{f:separable}
\end{figure}
The next subsections describe how the structure (\ref{e:wellsep}) can be used to construct a sparse factorization of $\mtx A$
that can be constructed and applied in $\mathcal O(N)$ time. 
First, Section \ref{s:skeleton_box} describes a choice of basis for (\ref{e:wellsep}), known as the skeleton basis,
where the set of source points in a box $\B$ is replaced by a subset that replicates the original effect in the far field of $\B$.
Such a basis can be computed for all boxes on a single level, leading to a sparse factorization of $\mtx A$, 
which is discussed in Section \ref{s:skeleton_level}.
Finally, Section \ref{s:multilevel_simple} discusses how this process can be applied recursively on a tree $\mathcal T$ to produce a multi-level algorithm
that involves only neighbor calculations. 
The key differences between our approach and traditional FMM algorithms are highlighted as well.

\subsection{Computing a skeleton basis for a single box} \label{s:skeleton_box}
Consider a box $\B_i$ with neighbors $\N$ and far-field $\F$. Abusing notations, we also use $\B_i$, $\N$, and $\F$ to denote the indices of points contained in the box, its neighbors, and its far field. Correspondingly, $\B_i\ \cup\ \N\ \cup\ \F$ forms a partitioning of row and column indices 
of the kernel matrix $\mtx A$.
For an appropriate permutation, we write $\mtx A$ as
\begin{equation}
    \label{ty_step1}
    \A
    = 
    \begin{pmatrix}
    \A_{\B_i \, \B_i} & \A_{\B_i \, \N} & \A_{\B_i \, \F} \\
    \A_{\mathcal{N} \, \mathcal B_i} & \A_{\mathcal{N} \, \mathcal{N}} & \A_{\mathcal{N} \, \mathcal F} \\
    \A_{\mathcal F \, \mathcal B_i} & \A_{\mathcal F \, \mathcal{N}} & \A_{\mathcal F \, \mathcal F}
    \end{pmatrix}.
\end{equation}
For non-oscillatory functions $G$, the matrices $\A_{\mathcal B_i \, \mathcal F}$ and $\A_{\mathcal F \, \mathcal B_i}$ for an arbitrary box $\B_i$ can be  approximated efficiently by a low-rank approximation for a prescribed accuracy $\eps$. In this work and other works based on recursive skeletonization, we use the interpolative decomposition (ID), which selects a subset of the original columns or rows of the matrix as a basis.
A (column) interpolative decomposition (ID) for a prescribed accuracy $\eps$ finds the so-called \emph{skeleton} indices $\S_i \subset \B_i$, the \emph{redundant} indices $\R_i = \B_i \backslash \S_i$, and an \emph{interpolation matrix} $\mtx{T}_{\B_i} \in \mathbb{C}^{|\S_i| \times |\R_i|}$ such that
\[
    \|\A_{\mathcal F \, \mathcal R_i} - \A_{\mathcal F \, \mathcal S_i} \, \mtx{T}_{\B_i}\| 
     \le \eps \, \|\A_{\mathcal F \, \mathcal B_i} \|.
\]
While the strong rank-revealing QR factorization of Gu and Eisenstat~\cite{gu1996efficient} is the most robust method for computing an ID, this work employs the column-pivoting QR factorization as a greedy approach~\cite{cheng2005compression}, which has better computational efficiency and  behaves well in practice. 
For a matrix of size $m \times n$ with numerical rank $k$, the cost to compute an ID using the aforementioned deterministic methods is $\bigO(mnk)$, which can be further reduced to $\bigO(mn \,\log k + k^2 n)$ using randomized sketching~\cite{dong2023simpler}.
Instead of compressing each block $\A_{\B_i \, \mathcal F}$ and $\A_{\mathcal F \, \B_i}$ separately, it is often convenient to conduct column ID compression of the concatenation, 
\begin{equation} \label{e:con}
\begin{pmatrix}
\A_{\mathcal F \, \mathcal B_i} \\
\A_{\mathcal B_i \, \mathcal F}^*
\end{pmatrix} 
= \begin{pmatrix}
\A_{\mathcal F \, \mathcal R_i} & \A_{\mathcal F\, \mathcal S_i} \\
\A_{\mathcal R_i \, \mathcal F}^* & \A_{\mathcal S_i\, \mathcal F}^* 
\end{pmatrix} \approx \begin{pmatrix}
 \A_{\mathcal F\, \mathcal S_i} \\
 \A_{\mathcal S_i\, \mathcal F}^* 
\end{pmatrix} \begin{pmatrix}
\mtx T_i & \mtx I
\end{pmatrix}
\end{equation}
which leads to a slightly larger set of skeleton indices but makes the implementation easier.
Notice that the computational cost would be $\bigO(N)$ if the full matrix in (\ref{e:con}) is formed, which turns out to be unnecessary. 
Instead, the effect of the far-field points can be replicated using an artificial surface of proxy points,
as is done in ~\cite{ying2004kernel,ying2006kernel,martinsson2005fast,ho2016hierarchical}; see 
\cite[Sec.~17.1]{martinsson2019fast} as well as \cite{ye2020analytical,xing2020interpolative} for a detailed discussion.
Instead of using the matrix in equation (\ref{e:con}), a smaller matrix $\begin{pmatrix}
\A_{{\rm proxy} \, \B_i} \\
\A_{\B_i \, {\rm proxy}}^*
\end{pmatrix}$ with $2\ n_{\rm proxy}$ rows is formed and compressed.

\begin{figure}
    \centering
    \includegraphics[width=0.45\textwidth]{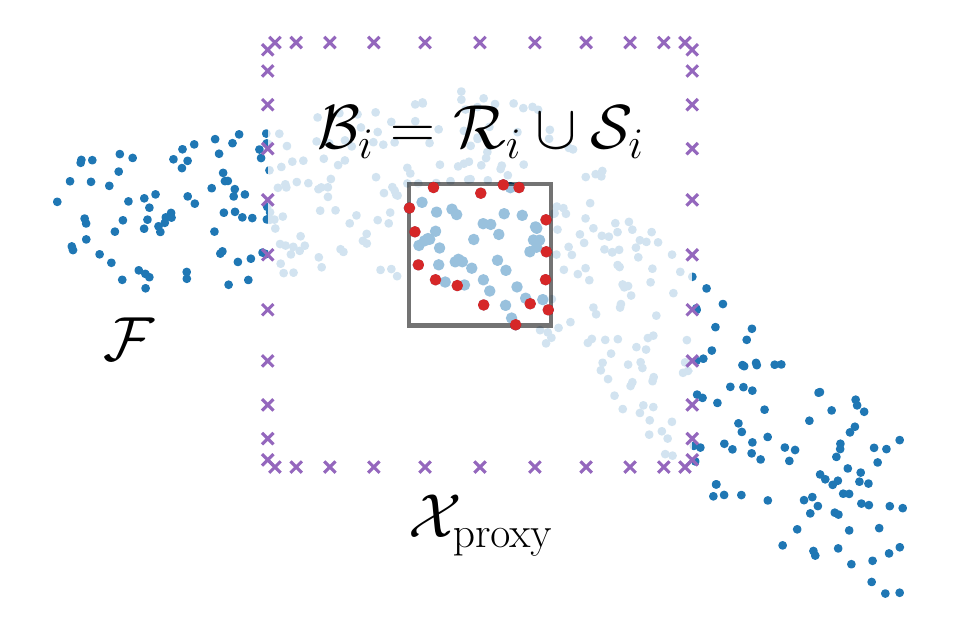}
    \caption{For a box $\B_i$, skeleton indices $\S_i \subseteq \B_i$ (shown in red) can be computed efficiently using
    a proxy surface, which are fictitious points replicating the effect of the
    far-field points ${\F}$. In our implementation of proxy surfaces, we use Chebyshev points on a box of length $2.95 a$ centered at $c$, where $a$ and $c$ are the length and center of box $\B_i$, respectively.}
    \label{f:proxy}
\end{figure}

With the appropriate choice of proxy surface, the matrix $\mtx A_{{\rm proxy}, \mathcal B_i}$ provides a low-dimensional representation that captures the essential features of the column space of $\mtx A_{\F,\B_i}$. Similarly, $\mtx A_{\B_i,{\rm proxy}}$ provides a low-dimensional representation of the row space of $\mtx A_{\B_i,\F}$.
While the two matrices have a different number of rows and columns, the proxy surface ensures that the essential relationship within these spaces is preserved for \textit{any} point distribution in the far-field. Therefore, the indices ${\B}_i={\R}_i \cup {\S}_i$ and interpolative matrix $\mtx T_i$ computed by forming and factorizing the smaller matrix will satisfy (\ref{e:con}).
Figure \ref{f:proxy} shows a proxy surface for a box $\B_i$ as well as the chosen subset of skeleton indices ${\S_i} \subseteq {\B_i}$.

\subsection{A Single-Level Factorization} \label{s:skeleton_level}
In this subsection, we use the interpolative decomposition to construct a sparse matrix
factorization for the finest level of a tree-based decomposition. The section also reviews existing 
approaches for kernel-independent fast multipole methods and contextualizes the simplified approach
proposed in this work.

Suppose that the matrix $\mtx A$ is tessellated into $p$ blocks $ \B_1, \dots \B_p$ of size $m = N/p$ according to the finest
level of $\mathcal T$.
The interpolative decomposition in equation (\ref{e:con}) can be computed for every box $\B_i$ using the techniques described
in Section \ref{s:skeleton_box}. 
Using (\ref{e:con}), the following equation holds for well-separated colleagues $\B_i, \B_j$ 
\begin{equation}\label{e:EFdecomp}
\mtx A_{\B_i,\B_j}  = 
\begin{pmatrix} \A_{\R_i, \R_j} & \A_{\R_i,\S_j}\\
\A_{\S_i, \R_j} & \A_{\S_i,\S_j}
\end{pmatrix}
\approx 
\begin{pmatrix}\mtx T_{\B_i}^*\\ \mtx I\end{pmatrix} 
\mtx A_{\S_i,\S_j}
\begin{pmatrix}\mtx T_{\B_j} & \mtx I\end{pmatrix} := 
\underset{m \times k}{\mtx E_{\B_i}}\ \underset{k \times k}{\mtx A_{\S_i,\S_j}}\ \underset{k \times m}{\mtx F_{\B_j}}.
\end{equation}
To define a single-level factorization of $\mtx A$, it is useful to separate $\mtx A$ into near and far interactions,
so that (\ref{e:EFdecomp}) can be used for the far interactions. We define $\mtx A^{\rm (near)}$ and $\mtx A^{\rm (far)}$ as
\begin{equation}
{\left[ \mtx A^{\rm (near)} \right]}_{i,j} := 
\begin{cases}\mtx A_{\B_i,\B_j}, &\text{if}\ \B_i, \B_j \ \text{are neighbors} \\
\mtx 0,\ &\text{otherwise}
\end{cases},\qquad \mtx A^{\rm(far)} := \mtx A - \mtx A^{\rm (near)}.
\end{equation}
It is also useful to define the set of all skeleton indices as 
\begin{equation}
{\S} = \bigcup_{i=1,\dots, p} {\S_i}, \qquad \text{where}\ \left| \S \right| = pk
\end{equation}
and corresponding matrices
$\underset{pk \times pk} {\mtx A_{\S,\S}}$, which is a submatrix of $\mtx A$, and 
$\underset{pk \times pk}{\mtx A_{\S, \S}^{\rm (near)}}$, where 
\begin{align}
{\left[ \mtx A_{\S,\S}^{\rm (near)} \right]}_{i,j} = 
\begin{cases}\mtx A_{\S_i,\S_j}, &\text{if}\ \B_i, \B_j \ \text{are neighbors}\\ \mtx 0, &\text{otherwise}. \end{cases}
\end{align}
Then $\mtx A^{\rm (far)}$ can be factorized as
\begin{equation}
\label{e:Atilde_formula}
\mtx A^{\rm (far)} \approx \left( \begin{matrix} \mtx E_{\B_1}\\& \ddots \\&& \mtx E_{\B_p}\end{matrix} \right)\ 
\left( \mtx A_{\S, \S} - \mtx A_{\S, \S}^{\rm (near)} \right)\ 
\left( \begin{matrix} \mtx F_{\B_1}\\& \ddots \\&& \mtx F_{\B_p}\end{matrix} \right) := \underset{N \times pk}{\mtx E}\ \underset{pk \times pk} {\widetilde {\mtx A}}\ \underset{pk \times N}{\mtx F}.
\end{equation}
The resulting factorization of $\mtx A$ is 
\begin{equation} \label{e:decomp_single}
\underset{N \times N}{\mtx A} \approx \underbrace{\underset{N \times N}{\mtx A^{\rm (near)}}}_{\text{(a) sparse}}\ +\ 
\underbrace{\underset{N \times pk}{\mtx E}\ \underset{pk \times pk}{\widetilde {\mtx A}}\ \underset{pk \times N}{\mtx F}}_{\text{(b) recursive term}}.
\end{equation}
The sparsity patterns of the matrices in equation (\ref{e:decomp_single}) are shown in Figure \ref{f:typical_fmm_step1} for the finest
level of the simple geometry in Figure \ref{f:1dgeom}.

\begin{figure}[htb!]
    \begin{subfigure}{\textwidth}
    \centering
    \input{FMM_typicalFMM_step1}
    \caption{Sparsity pattern of $\mtx A = \mtx A^{\rm (near)} + \mtx A^{\rm (far)}$, where $\mtx A^{\rm (far)}$ is factorized using operators $\mtx E$ and $\mtx F$. Note that $\widetilde {\mtx A}$ consists of sub-blocks of original kernel interactions of $\mtx A$.}
    \label{f:typical_fmm_step1}
    \end{subfigure}

    \vspace{0.5em}

    \begin{subfigure}{\textwidth}
    \centering
    \input{FMM_typicalFMM_step2prime}
    \caption{Sparsity pattern of the recursive term $\widetilde{\mtx A}$ which can be separated into the interaction list and the far interactions
    on the next coarser level.}
    \label{f:typical_fmm_step2}
    \end{subfigure}
    \caption{Previous approaches for the FMM separate the matrix $\mtx A$ into near and far interactions and factorize only
    the far-field interactions, as shown in Figure \ref{f:typical_fmm_step1}. 
    The approach can be applied recursively to $\widetilde{\mtx A}$, leading to the interaction list, shown in Figure \ref{f:typical_fmm_step2}.
    For the simple geometry of Figure \ref{f:1dgeom}, the interaction list for a box $\mathcal B$ is of size at most 3; generally, 
    interaction list is much larger for uniform point distributions in a square (at most size 27) or in a cube (at most size 189).}
    \label{f:fmm_typical}
\end{figure}

To generalize from a single-level decomposition (\ref{e:decomp_single}) 
to a multi-level algorithm,
a similar approach can be used recursively for $\widetilde{\mtx A}$.
The matrix can be tessellated into $4 \times 4$ blocks according to the
 tree decomposition on next coarser level.
However, the sparsity pattern of $\widetilde{\mtx A}$ is significantly different from the fully-dense matrix $\mtx A$.
Instead, $\widetilde{\mtx A}$ can be expressed as a sum of two matrices
\begin{equation}
\widetilde{\mtx A} = {\mtx A}_{\S,\S}^{\rm (interaction\ list)} + {\mtx A}_{\S,\S}^{\rm (far)},
\end{equation}
where the sparsity pattern is shown in Figure \ref{f:typical_fmm_step2}.
The matrix ${\mtx A}_{\S,\S}^{\rm (far)}$ consists of sub-blocks corresponding to well-separated boxes on the next coarser level.
The remaining matrix $\mtx A_{\S,\S}^{\rm (interaction\ list)}$ has the so-called interaction list.
Formally, box $\B_j$ is in the interaction list of box $\B_i$ when the boxes $\B_i$ and $\B_j$ are colleagues, the parents of $\B_i$ and $\B_j$ are neighbors, but the boxes are not themselves neighbors.
In 2D, the interaction list for a box $\B$ has size 27 for a fully populated tree,
and in 3D, it has size 189.
For the traditional FMM, additional data structures (e.g. W-list and X-list) are
needed for adaptive trees. These lists, like the interaction list, require additional bookkeeping and may cause difficulty in implementing the algorithm
efficiently.

\subsection{A Novel Approach using Modified Neighbor Translations}\label{s:multilevel_simple}

In this section, we present a modified approach which leads to a multi-level algorithm
which operates on much simpler data structures, requiring only the list of neighbor boxes at each level.
We focus on the case of a uniform multi-level tree (e.g. Figure \ref{f:1dgeom}) and discuss adaptive trees in Section \ref{s:algo};
however, adaptivity does not require much more machinery beyond what we discuss in this section.

Using the terminology of Section \ref{s:skeleton_level}, we can succinctly describe our modified approach by expanding the formula
(\ref{e:decomp_single}) using (\ref{e:Atilde_formula})
\begin{equation} \label{e:mod_near}
\underset{N \times N}{\mtx A^{\vphantom{\rm ()}}_{\vphantom{\S,}}}\
\approx \
\underbrace{ \underset{N \times N}{\mtx A^{(\rm near)}_{\vphantom{\S,}}}\
 - \
\underset{N \times pk}{\mtx E^{\vphantom{\rm ()}}_{\vphantom{\S,}}}\
\underset{pk \times pk}{\mtx A_{\S,\S}^{\rm (near)}}\
\underset{pk \times N}{\mtx F^{\vphantom{\rm ()}}_{\vphantom{\S,}}}\ }_{\text{(a) sparse}}
+\
\underbrace {\underset{N \times pk}{\mtx E^{\vphantom{\rm ()}}_{\vphantom{\S,}}}\
\underset{pk \times pk}{\mtx A^{\vphantom{\rm ()}}_{\S,\S}}\
\underset{pk \times N}{\mtx F^{\vphantom{\rm ()}}_{\vphantom{\S,}}}}_{\text{(b) recursive term}}.
\end{equation}
We make two observations regarding (\ref{e:mod_near}). First, the term (\ref{e:mod_near}a) is sparse and involves computations between neighbors
on the leaf level of the tree. Because the term $\mtx E\ \mtx A_{\S,\S}^{\rm (near)}\ \mtx F$ is explicitly subtracted in (\ref{e:mod_near}a),
the recursive term (\ref{e:mod_near}b) now has different structure, compared to the recursive term (\ref{e:decomp_single}b).
Second, in (\ref{e:mod_near}b), $\mtx A_{\S,\S}$ is the kernel matrix corresponding to all pairwise interactions among skeleton points. Therefore, $\mtx A_{\S,\S}$ has the same low-rank property as the original matrix $\mtx A$, and the same decomposition as in (\ref{e:mod_near}) can be applied recursively to $\mtx A_{\S,\S}$.

\begin{figure}[htb!]
    \begin{subfigure}{\textwidth}
    \centering
    \input{FMM_modifiedFMM_step1}
    \caption{In the modified approach, the neighbor calculations involve explicitly subtracting a sparse term. This leads to a fully dense matrix in the recursive term. The matrix $\mtx A_{\S,\S}$ is the kernel matrix corresponding to all pairwise interactions among skeleton points. Therefore, $\mtx A_{\S,\S}$ has the same low-rank property as the original matrix $\mtx A$ if it is viewed at a coarser level.}
    \label{f:modified_fmm_step1}
    \end{subfigure}

    \vspace{0.5em}

    \begin{subfigure}{\textwidth}
    \centering
    \input{FMM_modifiedFMM_step2}
    \caption{Sparsity pattern of the recursive term ${\mtx A}_{\S,\S}$. Note that it is fully dense, like the original matrix $\mtx A$.
    The matrix can be \textit{coarsened} and the same approach as in (\ref{e:mod_near}) can be used to factorize the recursive term.}
    \label{f:modified_fmm_step2}
    \end{subfigure}
    \caption{In the modified approach, we subtract an additional sparse term (\ref{e:mod_near}a), leading to a recursive term (\ref{e:mod_near}b)
    that has similar structure to the original matrix $\mtx A$.
    A multi-level algorithm is derived by coarsening the matrix 
    $\mtx A_{\S,\S}$ into 4 blocks, as shown in level 2
    of tree $\T$, of Figure \ref{f:1dgeom}. This leads to a multi-level
    algorithm involving \textit{modified} neighbor interactions at every level.}
    \label{f:fmm_modified}
\end{figure}

To describe the implementation of the matrix-vector multiply (\ref{e:matvec}) using the single-level decomposition (\ref{e:mod_near}), we introduce some notation using the language of potential theory.
With given charges $\mtx q$, we aim to compute the potential $\mtx u$, cf. Algorithm \ref{a:single_level}.
The terminology introduced is also used in Algorithm \ref{a:evaluate} for adaptive multi-level trees.

\begin{algorithm}[!htb]
\caption{\texttt{Apply using (\ref{e:mod_near})}}
\label{a:single_level}
  \begin{algorithmic}[1]
      \Require{\hspace{1em} Charges $\mtx q$.}
      \Ensure{Potentials $\mtx u$.}
      \State Compute $\hat{\mtx q}\quad :=\ \mtx F\ \mtx q$. 
      \Comment{\hfill \parbox[t]{.5\linewidth}{$\triangleright$ For box $\mathcal B_i$, the outgoing representation $\hat{\mtx q}_{\S_i} \in \mathbb{R}^k$ is an equivalent set of charges that approximates the effect of $\mtx q_{\B_i} \in \mathbb{R}^m$ in the far field. }}
      \State Set\quad \hspace{1.5em}\ ${\mtx q}_{\S}\ :=\ \hat {\mtx q}$.
      \State Compute $\mtx u_{\S}\ :=\ \mtx A_{\S,\S}\ {\mtx q}_{\S}$. 
      \Comment{\hfill \parbox[t]{.5\linewidth}{$\triangleright$ Recursively. }}
      \State Set\quad \hspace{1.5em}\ $\hat{\mtx u}\quad :=\ \mtx u_{\S}$.
      \State Update\quad $\hat {\mtx u}\quad \minuseq\ \mtx A_{\S,\S}^{\rm (near)}\ \hat {\mtx q}$. 
      \Comment{\hfill \parbox[t]{.5\linewidth}{$\triangleright$ For box $\mathcal B_i$, the incoming expansion $\hat{\mtx u}_{\S_i} \in \mathbb{R}^k$ approximately captures the effect of charges from the far field. }}
      \State Compute $\mtx u\quad :=\ \mtx A^{\rm (near)}\ \mtx q + \mtx E\ \hat{\mtx u}$. 
      \Comment{\hfill \parbox[t]{.5\linewidth}{$\triangleright$ Add contributions of neighbors.}}
  \end{algorithmic}
\end{algorithm}

The approach involves a slightly higher constant prefactor, compared to previous approaches, because (\ref{e:mod_near}) essentially involves
subtracting then adding a sparse matrix. The approach has benefit in ease of implementation as well as significantly less overhead to parallelizing operations
over the neighbor list instead of the interaction list (e.g. 7x fewer parallel tasks in 3D for the interaction list of size 189 and the neighbor list of size 27).
To implement the algorithm on a multi-level tree efficiently, 
we also require storing \textit{two}
types of incoming and outgoing expansions since algorithm involves subtracting
then summing terms; we elaborate further on these incoming and outgoing expansions
in Section \ref{s:algo}.

Briefly, we outline the scope of the remainder of the paper. In Section \ref{s:algo}, we describe the algorithm in detail for multi-level adaptive trees, as well as describe 
how the algorithm can be implemented efficiently in parallel. In Section \ref{s:complexity}, we analyze the complexity to precompute the skeleton points and apply the algorithm.
Finally, in section \ref{s:numerical}, we report numerical experiments on a variety of point distributions in 2D and 3D for the Laplace and low-frequency Helmholtz kernels.

\section{Algorithm}\label{s:algo}
This section describes the algorithm for adaptive quad-trees and oct-trees and extends the proposed methodology to non-uniform point distributions.
While Section \ref{s:methodology} describes a sparse factorization of $\mtx A$, this section
instead uses matrix-free notation (i.e.\ without explicitly forming any sparse matrices) which is more amenable to parallel implementation.
Additionally, the implementation of the precomputation stage and the FMM apply are discussed.

\subsection{Adaptive tree data structure} \label{s:tree}

The FMM relies on a hierarchical decomposition of the points $\mathcal{X} = \{x_i\}_{i=1}^N$ into an adaptive quad-tree or oct-tree.
Much of the terminology is the same as discussed earlier for the simple geometry in Figure \ref{f:1dgeom} which is defined on a uniform tree.
An adaptive tree can be formed using a recursive approach. First, the root box $\B_0$ which contains all the points is defined.
A parameter $b$ is also defined, for the maximum number of points a box may have.
Then $\B_0$ is split into 4 (or 8 for 3D) subdomains. Any box which is empty is pruned from the tree, and only boxes containing more than $b$ points
are subdivided recursively. The top-down construction produces a tree which is unbalanced. There may be adjacent leaf boxes that
vary tremendously in their size, leading to a possibly unbounded number of neighbors for a large box adjacent to very refined box.
To limit the number of interactions, adaptive trees are constrained to satisfy a 2:1 balance constraint, that is, adjacent leaf boxes must
be within one level of each other. Given an unbalanced tree, additional leaves may be added in a sequential procedure that produces a balanced tree;
as described in \cite{malhotra2015pvfmm}, this requires $\bigO(n_{\rm boxes} \log n_{\rm boxes})$, where $n_{\rm boxes}$ are the number of total boxes in the tree.

\begin{figure}[!htb]
\centering
\includegraphics[width=0.5\textwidth]{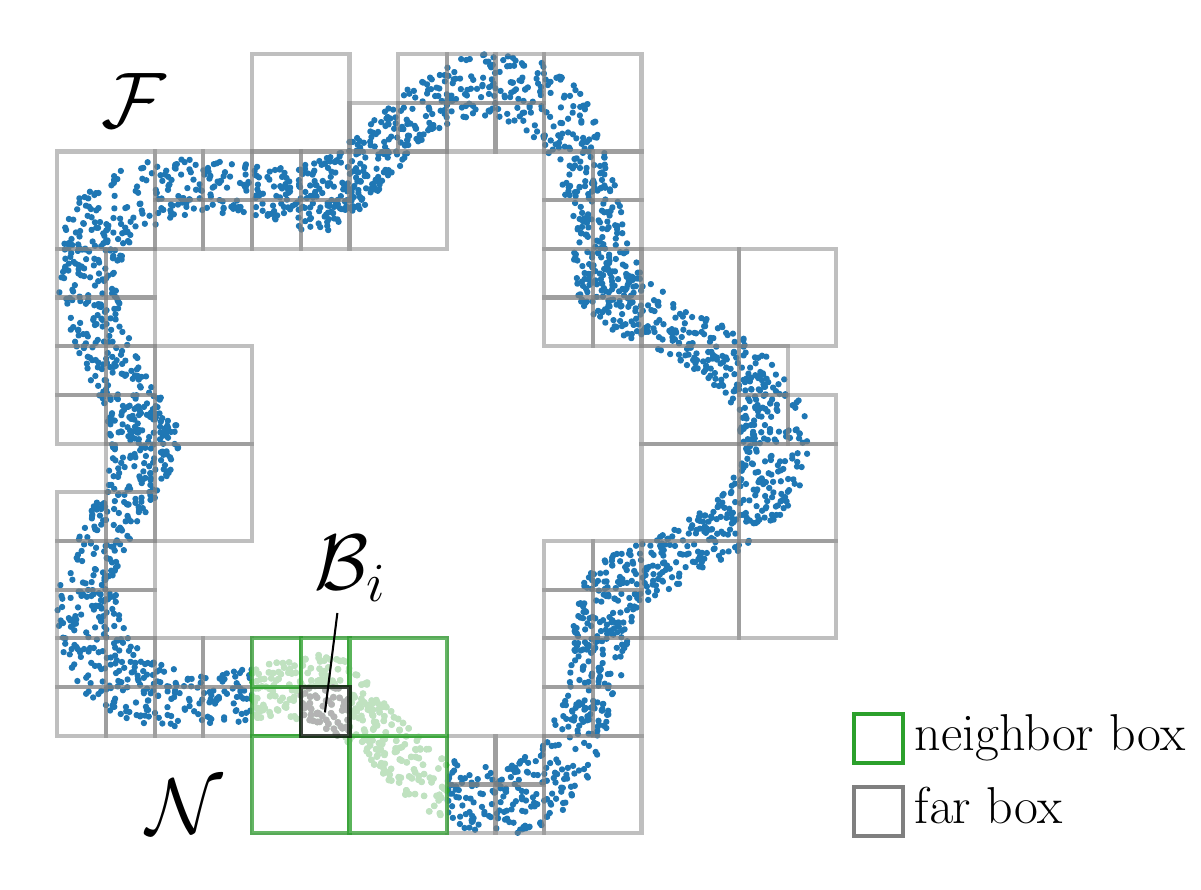}
\caption{The figure shows the leaf boxes of an adaptive quadtree for a non-uniform distribution of points.
The tree is adaptive, with reasonable restrictions on the adaptivity (e.g. adjacent leaf boxes are within one level of each other).
For box $\mathcal B_i$ with corresponding index set, we also show neighbor boxes in green and far boxes in gray.
The corresponding index sets for points in the near field and far-field are ${\N}$ and ${\F}$, respectively.
For adaptive trees, a box $\B_i$ may have a coarse neighbor on a level above or a fine neighbor
on a level below.}
\label{f:wiggly_geom}
\end{figure}

The \textit{near-field} $\mathcal N$ of box $\mathcal B$ of length $a$ and center $c$ contains points within an area defined by a box of length $3a$ centered at $c$.
For uniform trees, the near-field $\N$ only has adjacent colleagues of $\B$, that is, boxes on the same level as $\B$ which share a corner or an edge.
For adaptive trees, the near-field $\N$ may contain a \textit{coarse neighbor} on a level above
or a \textit{fine neighbor} on a level below. 
A coarse neighbor for a box $\B$ must necessarily be a leaf, which we denote $\L$;
for box $\B$ to have coarse neighbor $\L$, the parent $\P(\B)$ must be adjacent colleagues with $\L$, cf. Figure \ref{f:wiggly_geom}. Likewise, a leaf $\L$ may have a fine neighbor $\B$.
For $\T$ a uniform quad-tree or oct-tree, it is obvious that no box has more than 9 and 27 neighbor boxes, respectively.
For $\T$ an adaptive tree with a balance constraint, the number of neighbor boxes is bounded as well.
\subsection{Build stage}
As described in Section \ref{s:skeleton_box}, each leaf box $\B_i$ can be skeletonized into redundant and skeleton
indices $\B_i = \R_i \cup \S_i$ so that equation (\ref{e:con}) holds. The skeleton nodes $\S_i$ and interpolation 
matrix $\mtx T_{\B_i}$ that satisfies (\ref{e:con}) can be efficiently computed using a small proxy matrix.
After the skeleton nodes are computed for leaf boxes, skeleton indices are computed for the
remaining tree boxes (i.e.\ boxes with children) by traversing the tree upwards level by level.
For tree boxes $\B_i$, the indices $\B_i$ are set to be the union of skeleton indices of the children $\C(\B_i)$.
\input{alg_fact.tex}

\subsection{FMM apply}\label{s:fmm_apply}
Sections \ref{s:skeleton_level} and \ref{s:multilevel_simple} describe a sparse factorization of the kernel matrix $\mtx A$ in (\ref{e:matvec}) 
for the finest level of a uniform tree.
This section describes an efficient multi-level algorithm to compute (\ref{e:sum_fmm}) using a series of translations
between incoming and outgoing expansions for each box. The modifications needed for adaptive trees are also discussed.
Recall that some precomputation is required before computing (\ref{e:sum_fmm}); the skeleton indices for every box and interpolation matrices are
computed in Algorithm \ref{a:precompute} in an upward pass through the tree.

First, outgoing representations for the leaf boxes are computed as follows:
\begin{equation} \label{e:outgoing_leaf}
\hat{\mtx q}_{\S_i} := \begin{pmatrix}\mtx I &  \mtx T_{\L_i}\end{pmatrix}
\begin{pmatrix} {\mtx q}_{\S_i}\\ {\mtx q}_{\R_i} \end{pmatrix}, \qquad \text{where}\ \mtx A_{\F_i, \S_i} \hat {\mtx q}_{\S_i} \approx \mtx A_{\F,\B_i} \mtx q_{\B_i},\ \text{for all leaves}\ \L_i.
\end{equation}
In an upwards pass through the tree, outgoing representations $\mtx q_{\B_i}, \hat{\mtx q}_{\S_i}$ are computed for all tree boxes using formulas
\begin{align}
 \label{e:q_tree} \mtx q_{\B_i} &:= \left[ \hat{\mtx q}_{\S_j};  \hat{\mtx q}_{\S_{j+1}}; \dots \right]\ \text{for children}\ \B_j, \B_{j+1}, \dots\ \text{of box}\ \B_i\\
\label{e:qhat_tree} \hat{\mtx q}_{\S_i} &:= \begin{pmatrix}\mtx I &  \mtx T_{\B_i}\end{pmatrix}
\begin{pmatrix} {\mtx q}_{\S_i}\\ {\mtx q}_{\R_i} \end{pmatrix}.
\end{align}
The operation for computing ${\mtx q}_{\B_i}$ in (\ref{e:q_tree}) involves concatenating vectors, analogous to how the index vector ${\B_i}$ is accumulated from the children in Algorithm \ref{a:precompute} on line 6.

At level 1 of the tree, there are at most $2^d$ boxes, each of which contain $k\ 2^d$ active indices (assuming $k$ skeleton points per box), and
the potentials due to the skeleton charges can be computed directly
\begin{equation}
\mtx u_{\B_i}\ :=\ \sum_{\B_j\ \text{on level 1}} \mtx A_{\B_i, \B_j}\ \mtx q_{\B_j}, \qquad \text{for all boxes}\ \B_i\ \text{on level 1}.
\end{equation}

In a downward pass through the tree, outgoing expansions are computed for each box $\B_i$.
For each box $\B_i$ at level $l$, the following steps are performed.
First, initialize $\hat{\mtx u}_{\S_i}$ with the parent's contributions
    \begin{equation}
    \hat{\mtx u}_{\S_i}\ \text{$:= \mtx u_{\S_p}$}\ \text{from parent box $p = \P(\B_i)$,
    where $\mtx u_{\B_p}= \begin{pmatrix}\mtx u_{\S_p}\\ \mtx u_{\R_p}\end{pmatrix}$.}
    \label{e:init_hatu}
    \end{equation}
Then, subtract contributions of skeleton indices from colleague neighbors
    \begin{equation}
    \hat{\mtx u}_{\S_i}\ \minuseq\  \sum_{\substack{\text{colleague}\\ \text{neighbor}\ \B_j}} \mtx A_{\S_i,\S_j} \hat{\mtx q}_{\S_j}.
    \label{e:A_SS_near}
    \end{equation}
Finally, set the outgoing representation $\mtx u_{\B_i}$ to be the sum of contributions of box indices from colleague neighbors as well as the interpolated contribution of $\hat{\mtx u}_{\S_i}$
\begin{equation}
{\mtx u}_{\B_i} :=\ \begin{pmatrix}\mtx I \\ \mtx T_{\B_i}^*\end{pmatrix} \hat{\mtx u}_{\S_i}\ +\ 
\sum_{\substack{\text{colleague}\\ \text{neighbor}\ \B_j}} \mtx A_{\B_i,\B_j} \mtx q_{\B_j}.
\label{e:A_near}
\end{equation}
Equation (\ref{e:A_SS_near}) subtracts the contributions of the \textit{skeleton} indices of the neighbors
and the subsequent formula (\ref{e:A_near}) adds the contributions of the \textit{box} indices of the neighbors. 
There is a connection to these matrix-free formulas and the earlier discussion in Section \ref{s:methodology} that 
it is simple to see by referring to Figure \ref{f:modified_fmm_step1}.

For adaptive trees, a box may also have a fine or coarse neighbor on the level below or above, respectively.
This requires modifying the formulas (\ref{e:A_SS_near}) and (\ref{e:A_near}).
We again proceed level-by-level. For all boxes on a level, we initialize $\hat{\mtx u}_{\S_i}$ as described
in equation (\ref{e:init_hatu}).
Then, we subtract contributions from coarse neighbors, as well as colleague neighbors
\begin{equation}
    \hat{\mtx u}_{\S_i}\ \minuseq\ \left( \sum_{\substack{\text{colleague}\\ \text{neighbor}\ \B_j}} \mtx A_{\S_i,\S_j} \hat{\mtx q}_{\S_j}
    +\ \sum_{\substack{\text{coarse}\\ \text{neighbor}\ \B_j}} \mtx A_{\S_i, \B_j} \mtx q_{\B_j} \right).
\end{equation}
Finally, in addition to the terms in (\ref{e:A_near}), there are contributions from fine and coarse neighbors
\begin{align}
\begin{split}
{\mtx u}_{\B_i} &:=\quad \begin{pmatrix}\mtx I \\ \mtx T_{\B_i}^*\end{pmatrix} \hat{\mtx u}_{\S_i}\ +\ 
\sum_{\substack{\text{colleague}\\ \text{neighbor}\ \B_j}} \mtx A_{\B_i,\B_j} \mtx q_{\B_j}\\\
 &+\ \sum_{\substack{\text{coarse}\\ \text{neighbor}\ \B_j}} \mtx A_{\B_i, \B_j} \mtx q_{\B_j} +
\sum_{\substack{\text{fine}\\ \text{neighbor}\ \B_j}} \left( \mtx A_{\B_i, \B_j} {\mtx q}_{\B_j} -  \mtx A_{\B_i, \S_j} \hat{\mtx q}_{\S_j} \right).
\end{split}
\end{align}
The resulting potential $\mtx u$ is obtained by assembling the potentials $\mtx u_{\L}$ for every leaf box $\L$ into a single long vector. 

\subsection{Parallel Implementation}\label{s:gpu_implementation}

To describe the efficient parallel implementation of the algorithm, we reorder the computation and introduce new notations. 
In particular, we define translation operators that simplify the algorithm's description. 
The key idea is that the order of operations can be rearranged so that the algorithm proceeds as follows:
\begin{itemize}
    \item  Upward pass. Compute the outgoing expansions for each box using the precomputed interpolation matrices $\mtx T_{\B_i}$. These matrices transfer the outgoing expansion from child boxes to their parent boxes.
    \item Parallel translations between the neighbor list. Perform an embarrassingly parallel computation to determine the incoming expansion for each box. This step uses the incoming-from-outgoing translation operators between each box and its colleague neighbors. For a non-uniform tree, additional operations are required, including incoming-from-sources and targets-from-outgoing translations for both coarse and fine neighbor boxes.
    \item Downward pass. Compute the incoming expansions for each box. Once again, the interpolation matrices $\mtx{T}_{\B_i}$ are used to transfer incoming expansions from parent boxes to their children, ultimately yielding the computed potential.
\end{itemize}
Reordering the computation in this manner ensures that the only {level-to-level dependence} occurs in the upward and downward passes. 
The translation operations, which involve evaluating kernel matrix entries, are the most computationally intensive portion of the algorithm. 
However, their {embarrassingly parallel nature} makes them well-suited for parallelization.

Table~\ref{tab:notation_apply} summarizes the notation used to describe {Algorithm~\ref{a:evaluate}}. 
The parallel translation stage relies on the use of \textit{translation operators}. 
For example, the incoming-from-outgoing translation converts outgoing representations 
\((\mtx{q}_{\B_j}, \hat{\mtx{q}}_{\S_j})\) in box \(\B_j\) to incoming representations 
\((\mtx{u}_{\B_i}, \hat{\mtx{u}}_{\S_i})\) in box \(\B_i\), where boxes \(\B_i\) and \(\B_j\) are colleague neighbors.

Specifically, for a box \(\B_i\) with neighboring colleague \(\B_j\), the operator is defined as \textit{adding} to the outgoing expansion \(\mtx{u}_{\B_i}\) and \textit{subtracting} from the outgoing expansion \(\hat{\mtx{u}}_{\S_i}\) by reordering equations (21a) and (22a). The operation is expressed as:
\begin{equation*}
    (\mtx{u}_{\B_i}, \hat{\mtx{u}}_{\S_i}) 
    \xleftarrow{\rm (ifo)} (\mtx{q}_{\B_j}, \hat{\mtx{q}}_{\S_j}),
    \qquad
    \mtx{u}_{\B_i} \pluseq \mtx{A}_{\B_i, \B_j} \mtx{q}_{\B_j}, 
    \quad
    \hat{\mtx{u}}_{\S_i} \minuseq \mtx{A}_{\S_i, \S_j} \hat{\mtx{q}}_{\S_j}.
\end{equation*}
For non-uniform trees, additional translation operators are introduced. The translation operators
$\B_i \xleftarrow{\rm (ifs)} \L_j$, and 
$\L_i \xleftarrow{\rm (tfo)} \B_j$
are used to translate expansions between a box with a fine or coarse neighbor.
These translation operators are defined in {Algorithm~\ref{a:evaluate}}.

\begin{table}[!htb]
    \centering
    \begin{tabular}{|c|l|}
    \hline
         $\mtx q_{\L_i}$ & charges at sources $\X_{\L_i}$\\
         & \hspace{1em} (given at leaves)\\ \hline
         $\mtx u_{\L_i}$ & potentials at targets $\X_{\L_i}$\\
         & \hspace{1em} (algorithm output)\\ \hline
         $\mtx q_{\B_i}, \hat{\mtx q}_{\S_i}$ & outgoing representations for box $\B_i$\\
         $\mtx u_{\B_i}, \hat{\mtx u}_{\S_i}$ & incoming representations for box $\B_i$\\ \hline
         $\B_i \xleftarrow{\rm (ifo)} \B_j$& translation operator $({\mtx u}_{\B_i}, \hat{\mtx u}_{\S_i}) \leftarrow (\mtx q_{\B_j}, \hat {\mtx q}_{\S_j})$\\
         & \hspace{1em} $\mathcal B_i \leftarrow \text{colleague neighbor}\ \mathcal B_j$\\
         $\B_i \xleftarrow{\rm (ifs)} \L_j$ & translation operator $(\mtx u_{\B_i},\hat{\mtx u}_{\B_i}) \leftarrow {\mtx q}_{\L_j}$\\
         & \hspace{1em} $\mathcal B_i \leftarrow \text{coarse neighbor}\ \L_j$\\
         $\L_i \xleftarrow{\rm (tfo)} \B_j$ & translation operator ${\mtx u}_{\L_i} \leftarrow ({\mtx q}_{\B_j}, \hat {\mtx q}_{\S_j})$\\
         & \hspace{1em} $\L_i \leftarrow \text{fine neighbor}\ \B_j$\\
         \hline
    \end{tabular}
    \caption{Notation used to describe the FMM apply. $\B_i$ is a box in the tree (may be a leaf or not), and $\L_i$ is specifically a leaf box.}
    \label{tab:notation_apply}
\end{table}

\input{alg_solve.tex}

\section{Complexity analysis} \label{s:complexity}
In this section, the computational complexity and storage costs of 
Algorithms \ref{a:precompute} and \ref{a:evaluate} are analyzed.
Assume that the point distribution of $N$ particles is
partitioned into a tree with at most $b$ points in the leaf boxes, so that the total number of boxes is
\begin{equation}
n_{\rm boxes} := n_{\rm leaf} + n_{\rm tree}, \qquad n_{\rm leaf} \approx N/b,\quad  n_{\rm tree} \approx \frac{1}{2^d-1} \frac N b.
\end{equation}
The number of neighbor boxes $n_{\rm neigh}$ 
depends on the point geometry; for a uniform point distribution, $n_{\rm neigh}=3^d$, and the quantity is bounded for balanced adaptive trees.
We also define
\begin{equation}
\mathcal N_{\rm coarse}(\B_i) := \{ \B_j\ |\ \P(\B_i)\ \text{and}\ \P(\B_j)\ \text{are neighbors} \}, \quad 
n_{\rm coarse} := \left| {\N_{\rm coarse}(B_i)} \right|
\end{equation}
as the set of boxes such that the \textit{parents} are neighbors. Note that this list is a superset of the interaction list
and can also be defined on the next coarse level of the tree.
The size of this list depends on the point distribution;
it is of size at most $n_{\rm coarse} = 6^d$ for uniform trees.
We also define constants $t_{\rm kernel}$ and $t_{\rm flop}$
as the time needed for one kernel evaluation and one linear algebraic flop, respectively. 
We also assume that the numerical rank from ID compression is bounded by a constant $k_{\rm max}$.

\vspace{0.5em}

\noindent \textbf{Precomputation Costs}:
For each box, we evaluate the kernel interactions between indices $\B_i$ and proxy surface
then compute the column ID.
Assuming the accumulated skeleton indices have size at most $b$ and the proxy surface has at 
most $n_{\rm proxy}$ points, then the total cost is
\begin{align*}
    T_{\rm skel} &= n_{\rm leaf}\ n_{\rm proxy}\ \bigO \left( {b\ t_{\rm kernel}\ + b\ k_{\rm max}\ t_{\rm flop}} \right)\\
    &+\ n_{\rm tree}\ n_{\rm proxy}\ \bigO \left( {2^d\ k_{\rm max}\ t_{\rm kernel} + 2^d\  k_{\rm max}^2\ t_{\rm flop}} \right) \\
    &= \bigO \left(  n_{\rm proxy}\ (t_{\rm kernel} +  k_{\rm max}\ t_{\rm flop})\ N \right), 
\end{align*}
assuming that $b \ge k_{\rm max}$. Choosing the leaf size to be at least $k_{\rm max}$ is a reasonable guideline to achieve
competitive memory complexity as well. Under this assumption,
the memory needed to store $\mtx T_{\B_i}$ for all boxes is
\begin{align*}
    M_{\rm proj} &= \bigO \left(  k_{\rm max}\ N \right).
\end{align*}

\noindent \textbf{Algorithm Apply}:
The first stage of the algorithm requires an upward pass to compute outgoing expansions $\mtx q, \hat{\mtx q}$ (referred to $T_{\rm ofo}$
for the time to translate outgoing expansions).
Then, computing incoming  expansions $\mtx u, \hat{\mtx u}$ requires kernel evaluations between every box and its list of neighbors (referred to as
$T_{\rm ifo}$ for the time for translate outgoing to incoming expansions).
Finally, there is a downward pass to compute the potential at the leaves (referred to as $T_{\rm ifi}$ for the time to translate incoming expansions). 
The total time needed is
\begin{equation}
T_{\rm apply} = T_{\rm ofo} + T_{\rm ifo} + T_{\rm ifi}.
\end{equation}
Both the upward and downward passes require matrix-vector
products with the pre-computed interpolation operators $\mtx T_{\B_i}$, so 
\begin{equation}
T_{\rm ofo} + T_{\rm ifi} = \bigO \left( k_{\rm max}\ t_{\rm flop}\ N \right).
\end{equation}
We now analyze the costs for translations $T_{\rm ifo}$.
At the leaf level, the interactions between neighbors are computed directly; this requires
$n_{\rm leaf}\ n_{\rm neigh}\ b^2$ kernel evaluations as well as linear algebraic flops.
The incoming expansions $\mtx u, \hat{\mtx u}$ for each box (i.e. tree or leaf box) is computed by first subtracting interactions 
over the skeleton indices of the neighbor list $\N(\B_i)$,
then adding interactions over the skeleton indices of $\N_{\rm coarse}(\B_i)$. The results of these computations
are accumulated into $\mtx u, \hat{\mtx u}$ and added in the downward pass, as described in Algorithm \ref{a:evaluate}. 
A key observation of this work is that the traversal over $\N_{\rm coarse}(\B_i)$
can be regrouped and described as a traversal of the neighbors \textit{at the next coarse level}.
When implementing the subtraction of neighbor interactions, one can reuse matrices that have already been formed, saving on kernel evaluation costs.
The total cost is then
\begin{align}
\begin{split}
    T_{\rm ifo} &= n_{\rm leaf}\ n_{\rm neigh}\ b^2\ (t_{\rm kernel} + t_{\rm flop}) \\
    &+ n_{\rm boxes}\ \left( n_{\rm coarse}\ k_{\rm max}^2\ (t_{\rm kernel} + t_{\rm flop})\ +\ n_{\rm neigh}\ k_{\rm max}^2\ t_{\rm flop} \right)\\
    T_{\rm apply} &= \bigO \left(\ \left({b\ n_{\rm neigh} + k_{\rm max}\ n_{\rm coarse}} \right)\ \left({t_{\rm kernel} + t_{\rm flop}} \right)\ N\ \right).
\end{split}
\end{align}
The resulting asymptotic costs are slightly higher than the traditional FMM, however, by grouping $\N_{\rm coarse}(\B_i)$ on the coarser level of the tree,
the algorithm traverses the neighbor list at every level. This leads to fewer parallel tasks (e.g. smaller overhead cost) and substantial 
benefits in the ease of implementing the algorithm.

\section{Numerical results} \label{s:numerical}

This section demonstrates the performance of a parallel implementation of the algorithm for both {uniform} and {non-uniform} point distributions in {2D} and {3D}. The implementation consists of a {precomputation stage} on the CPU. Once precomputation is complete, the algorithm rapidly computes the potential (\ref{e:sum_fmm}) on the GPU.

\subsection{Precomputation Stage}

The precomputation stage determines skeleton indices and computes interpolation
matrices $\mtx T_{\B_i}$, which translate incoming and outgoing expansions between child and parent boxes. This phase is performed on the CPU in double precision (or complex double precision for Helmholtz kernels), using OpenMP for parallelization across boxes at each tree level. The skeleton rank  $k$
adaptively chosen for each box based on a user-specified tolerance, with the maximum rank across all boxes denoted as $k_{\rm max}$. 

Because the apply-stage uses single precision, these are stored in single precision (or complex single precision for Helmholtz kernels). We report precomputation metrics $T_{\rm skel}$ as the time required to compute the skeleton indices and interpolation matrices, as well 
as $M_{\rm proj}$, which is the memory needed to store all interpolation matrices $\mtx T_{\B_i}$. 

The method operates on an {adaptive tree} $\mathcal{T}$ that satisfies the {2:1 balance constraint}, as described in Section \ref{s:tree}. Tree construction begins with a downward traversal to create an unbalanced tree where each leaf box contains at most $b$ points. The tree is then balanced by introducing additional leaf nodes using a sequential algorithm. Geometric lookups are handled efficiently using Morton codes. The time to partition and balance the tree is denoted as $T_{\rm tree}$.

\subsection{GPU Implementation of Apply Stage}

Once the precomputation is complete, the apply stage computes the potential (\ref{e:sum_fmm}) using the GPU for computationally intensive tasks. The computation is conducted in single precision (or complex single precision for the Helmholtz kernels).

The upward pass, which computes outgoing expansions, relies on matrix-vector products
with the precomputed interpolation matrices $\mtx T_{\B_i}$. To reduce
memory overhead, the interpolation matrices are stored on the CPU, and the upward pass is executed on the CPU. The computed outgoing expansions are then transferred to the GPU.

The {translation operators}, which convert outgoing expansions to incoming expansions, are the most computationally intensive portion of the FMM apply stage. These operations involve kernel evaluations and matrix-vector products between sub-blocks of the kernel matrices $\mtx A_{\B_i,\B_j}$. To optimize GPU performance, we use the following strategies:
\begin{itemize}
    \item Batching. Translation operators are grouped into large batches of tens of thousands of submatrices to exploit GPU parallelism efficiently. 
    \item Zero padding. Adaptive ranks are padded with zeros to ensure uniform vector lengths, enabling efficient batched operations.
    \item On-the-fly row generation. Instead of forming the submatrices $\mtx A_{\B_i,\B_j}$ explicitly (which incurs significant memory overheads), rows are generated on-the-fly in batches and applied to the vector sequentially, row by row. This approach minimizes memory overheads and is highly cache-optimized, especially when processing over large batch sizes of grouped submatrices.
\end{itemize}

Once the outgoing expansions are computed, they are transferred back to the CPU for the downward pass, which propagates incoming expansions from the parents to the children using the interpolation matrices. The result is accumulated values $\mtx u_{\L_i}$ stored in the leaf boxes. The time needed to compute the potential using Algorithm \ref{a:evaluate} is reported as $T_{\rm apply}$, as well as \text{relerr}, the maximum relative error, evaluated on a subset of points. Table \ref{t:notation} summarizes the notation used for reporting the results.

\begin{table}[!htb]
    \centering
    \begin{tabular}{|c|l|} \hline
     $N$   & number of points \\
     $b$   & leaf size of tree\\
     $k_{\rm max}$   & maximum skeleton rank \\ \hline
     $T_{\rm tree}$ & time to partition and balance tree\\
     $T_{\rm skel}$ & time for Algorithm \ref{a:precompute}\\
     $M_{\rm proj}$ & memory needed for $\mtx T_{\mathcal B}$ for all boxes\\
     $T_{\rm apply}$ & time for Algorithm \ref{a:evaluate}\\ \hline
     \multirow{2}{*}{relerr} & maximum relative error\\
     & (evaluated on a subset of points)\\ \hline
    \end{tabular}
    \caption{Notation for reported numerical results.}
    \label{t:notation}
\end{table}

\subsection{Overview of Implementation}

The code is implemented in \texttt{Python}. The pre-processing stage uses the packages \texttt{numpy} and \texttt{scipy} with parallelization via OpenMP on the CPU.

The GPU implementation of the FMM apply stage uses \texttt{PyTorch}. Most of the computation time is spent on evaluating and applying kernel matrices. To optimize these operations, the \texttt{Pykeops} package~\cite{JMLR:v22:20-275} is employed. \texttt{Pykeops} generates highly optimized, compiled CUDA implementations of batched matrix-vector multiplications with kernel matrices, ensuring minimal memory overhead.
The Hankel function \(H_0^{(1)}\) is not available in \texttt{Pykeops}, and the translation operators for the 2D Helmholtz kernel are implemented (less efficiently) in \texttt{Cupy}.

The numerical experiments were conducted on a workstation equipped with an Intel Xeon Gold 6326 CPU (2.9 GHz, 16 cores) with 250~GB of RAM as well as an NVIDIA A100 GPU with 80~GB of RAM. The code is publicly available at
\cite{Yesypenko_SkelFMM_A_Simplified_2025}
with tutorials for implementing the algorithm.

\subsection{2D Experiments}
Figure~\ref{f:2dlp} presents results for the Laplace kernel, while Figure~\ref{f:2dhh} reports results for the 2D Helmholtz kernel with \(\kappa = 100\). Both figures illustrate the performance of the method across a range of user-specified tolerances. For 2D point distributions, the skeleton rank scales approximately as \(k_{\rm max} \sim \log(1 / \epsilon)\) for the Laplace equation.

The method shows strong performance on the GPU, with the runtime \(T_{\rm apply}\) increasing only slightly as the user tolerance \(\epsilon\) is refined. However, for small sub-matrix sizes, there may be overheads associated with GPU usage. The runtime remains nearly constant in the pre-asymptotic range before transitioning to a linear scaling regime.

\input{results_2d_lp}

\input{results_2d_hh}

\input{results_3d_lp}
\input{results_3d_hh}

\input{results_comparison}

\subsection{3D Experiments}
Figures~\ref{f:3dlp} and~\ref{f:3dhh} present results for the Laplace and Helmholtz kernels \(\kappa = 20\), respectively, for 3D point distributions. For points distributed within 3D volumes, the skeleton ranks \(k_{\rm max}\) increase significantly as a function of the tolerance \(\epsilon\), approximately scaling as $k \sim \left( \log (1/ \epsilon) \right)^2$. The oct-tree data structure also 
leads to larger accumulated tree box sizes, increasing from 
$4k$ in 2D to $8k$ in 3D. Despite the increased computational demands due to higher ranks, larger sub-block sizes, and a greater number of near neighbors in 3D, the runtime for the 3D Laplace equation remains competitive. This efficiency is largely due to the effective GPU handling of translation operators. Additionally, the larger matrix sub-block sizes in 3D benefit from GPU parallelism, mitigating the overheads typically observed with smaller sub-matrix sizes in 2D.

Selecting skeleton points as a subset of the original points is particularly effective for representing surface distributions. As shown in Table~\ref{t:sphere_lp}, this method achieves high accuracy with relatively low skeleton ranks, especially when compared to volumetric distributions (Table~\ref{t:cube_lp}). The cost of computing a skeleton representation tailored to the geometry can be amortized across multiple right-hand sides, making it particularly useful for the iterative solution of discretized boundary integral equations.

\subsection{Performance Comparison}

We compare the performance of \texttt{FMM3D}, a widely used CPU implementation of the fast multipole method, to \texttt{skelFMM} in its CPU and GPU variants (Figure~\ref{f:3dlp_comp}). Table~\ref{t:3dlp_comp} provides a detailed breakdown of the time spent in various parts of the \texttt{skelFMM} algorithm. \texttt{FMM3D} \cite{fmm3d} employs analytic expansions and advanced techniques such as directional groupings and diagonal translations~\cite{greengard1997new, rokhlin1992diagonal}, making it a good benchmark for assessing the performance of alternative FMM implementations.

For uniform distributions in the volume, the results indicate room for improving \texttt{skelFMM}'s performance, as \texttt{FMM3D} incorporates
advanced techniques, such as directional groupings, to more efficiently process the interaction list.
In contrast, for adaptive distributions, such as on the sphere surface, the performance gap between the CPU implementations of the two codes much narrower, indicating effective handling of non-uniform point distributions by \texttt{skelFMM} and benefits from \texttt{skelFMM}'s tailored skeleton representation.

The results clearly demonstrate the substantial performance improvements achieved through GPU acceleration. On the CPU, processing large subblock sizes incurs significant computational overhead, whereas the GPU efficiently manages many thousands of subblocks using batched operations.
Notably, computational bottlenecks present in the CPU implementation of \texttt{skelFMM} are largely mitigated on the GPU. \texttt{skelFMM cuda} achieves significant speedups, outperforming \texttt{skelFMM cpu} by 34x for the cube and 20x for the sphere. \texttt{skelFMM cuda} also outperforms \texttt{FMM3D} by 11x for the cube and 40x for the sphere, though these results should be interpreted with caution, as they reflect the differences in CPU and GPU architectures rather than a purely algorithmic advantage.

\section{Conclusions}

This work introduces a novel kernel-independent algorithm for the fast multipole method (FMM). Similar to many other kernel-independent methods, analytical expansions are replaced by a chosen set of `skeleton' points, which replicate the effect of the original source points in the far field. 

A key contribution of this work is the elimination of interaction lists and auxiliary data structures from the kernel-independent FMM, while retaining compatibility with adaptive trees. Instead, the proposed algorithm describes novel incoming-from-outgoing translation operators involving only the neighbor list at each level, resulting in much simpler data structures. 
Although this adjustment introduces a slight increase in computational complexity, it significantly streamlines the computation, reducing parallel overheads and enhancing the overall efficiency of the algorithm.

Numerical results confirm the effectiveness of the GPU implementation for the Laplace and Helmholtz kernels across a range of tolerances and problem sizes. The GPU-accelerated \texttt{skelFMM cuda} variant delivers substantial speedups compared to its CPU counterpart and exhibits competitive performance with benchmark code \texttt{FMM3D}, particularly for adaptive point distributions.

This framework opens exciting opportunities for further development. For volume distributions or scenarios requiring faster precomputation, adopting a uniform skeleton basis across boxes, as explored in~\cite{malhotra2015pvfmm, malhotra2016algorithm}, could further enhance efficiency. FFT-based accelerations and directional groupings could also optimize translation operations, particularly for high-accuracy 3D computations. Additionally, the simplicity of the neighbor-list-based approach positions this method as a strong candidate for extensions into efficient direct solvers or preconditioners for integral equations, where simpler data structures can lead to significant computational benefits.

\vspace{1em}

\subsection*{Acknowledgments}
The work reported was supported by the Office of Naval Research (N00014-18-1-2354),
by the National Science Foundation (DMS-1952735, DMS-2012606, and DMS-2313434),
and by the Department of Energy ASCR (DE-SC0022251).
We thank Umberto Villa for access to computing resources.

%% file: 1d_geom.tex
\begin{tikzpicture}[scale=0.7]

\pgfmathsetmacro{\b}{1}
\pgfmathsetmacro{\m}{8}

\pgfmathsetmacro{\offsetx}{0}
\pgfmathsetmacro{\offsety}{0}

\filldraw[fill=none, thick] (\offsetx, \offsety+\m*\b) rectangle (\offsetx + \m*\b, \offsety+\m*\b + \b);
\draw plot [only marks, mark=*, mark size=0.5, domain=0:8, samples=700] ({\x} ,{rnd + 0.5 + 7.5});
\node at (\offsetx-2*\b, \offsety+\m+0.5*\b) {Points $\{x_i\}_{i=1}^N$};

\pgfmathsetmacro{\offsetx}{0}
\pgfmathsetmacro{\offsety}{-1.5*\b}
\pgfmathsetmacro{\mtmp}{\m* 0.25}
\node at (\offsetx-1.5*\b, \offsety+\m+0.5*\b) {Level 1};
\filldraw[fill=none, thick] (\offsetx, \offsety+\m*\b) rectangle (\offsetx + \m*\b, \offsety+\m*\b + \b);
\foreach \x in {1,...,\mtmp}{

\pgfmathtruncatemacro{\xtmp}{\x+1}
\draw[thick] (\offsetx + 4*\x*\b, \offsety+\m*\b) -- (\offsetx + 4*\x*\b , \offsety+\m*\b+\b);
\node[scale=0.8] at (\offsetx + 4*\x*\b - 2*\b, \offsety+\m*\b+0.5*\b) {${\xtmp}$};
}

\pgfmathsetmacro{\offsetx}{0}
\pgfmathsetmacro{\offsety}{-3*\b}
\pgfmathsetmacro{\mtmp}{\m* 0.5}
\node at (\offsetx-1.5*\b, \offsety+\m+0.5*\b) {Level 2};
\filldraw[fill=none, thick] (\offsetx, \offsety+\m*\b) rectangle (\offsetx + \m*\b, \offsety+\m*\b + \b);
\foreach \x in {1,...,\mtmp}{

\draw[thick] (\offsetx + 2*\x*\b, \offsety+\m*\b) -- (\offsetx + 2*\x*\b , \m*\b + \offsety+\b);
\pgfmathtruncatemacro{\xtmp}{\x+3}
\node[scale=0.8] at (\offsetx + 2*\x*\b - \b, \offsety+\m*\b+0.5*\b) {${\xtmp}$};
}

\pgfmathsetmacro{\offsetx}{0}
\pgfmathsetmacro{\offsety}{-4.5*\b}
\node at (\offsetx-1.5*\b, \offsety+\m+0.5*\b) {Level 3};
\filldraw[fill=none, thick] (\offsetx, \offsety+\m*\b) rectangle (\offsetx + \m*\b, \offsety+\m*\b + \b);
\foreach \x in {1,...,\m}{
\draw[thick] (\offsetx + \x * \b, \offsety+\m*\b) -- (\offsetx + \x*\b , \offsety+\m*\b + \b);
\pgfmathtruncatemacro{\xtmp}{\x+7}
\node[scale=0.8] at (\offsetx + \x*\b - 0.5*\b, \offsety+\m*\b+0.5*\b) {${\xtmp}$};
}
\end{tikzpicture}

%% file: FMM_separable.tex
\begin{tikzpicture}[scale=0.14]

\pgfmathsetmacro{\b}{2}
\pgfmathsetmacro{\k}{1}

\pgfmathsetmacro{\m}{8}

\pgfmathsetmacro{\offsetx}{2.5*\m}
\pgfmathsetmacro{\offsety}{0}

\filldraw[fill=gray, thick] (\offsetx,\b) rectangle (\offsetx+\m*\b,\m*\b+\b);
\foreach \x in {1,...,\m}
\draw[thick,white,fill=darkgray] (\offsetx + \m*\b - \x * \b + \b, \offsety + \x * \b + \b) rectangle (\offsetx + \m*\b - \x * \b , \offsety + \x * \b );
\foreach \x in {2,...,\m}
\draw[thick,white,fill=darkgray] (\offsetx + \m*\b - \x * \b + 2*\b,\offsety+ \x * \b + \b) rectangle (\offsetx+\m*\b - \x * \b + \b, \offsety+\x * \b );
\foreach \x in {2,...,\m}
\draw[thick,white,fill=darkgray] (\offsetx+\m*\b - \x * \b + \b, \offsety+\x * \b) rectangle (\offsetx+\m*\b - \x * \b, \offsety+\x * \b - \b);

\foreach \x in {1,...,\m}
\draw[thick, white] (\offsetx,\offsety+\x*\b) -- (\offsetx+2*\m,\offsety+\x*\b); 
\filldraw[fill=none, thick] (\offsetx,\offsety+\b) rectangle (\offsetx+\m*\b,\offsety+\m*\b+\b);

\pgfmathsetmacro{\offsetx}{5*\m}
\pgfmathsetmacro{\offsety}{0}

\filldraw[fill=gray, thick] (\offsetx +0,\offsety+\b) rectangle (\offsetx +\m*\b,\offsety+\m*\b+\b);

\foreach \x in {1,...,\m}
\draw[thick, white,fill=darkgray] (\offsetx + \m*\b - \x * \b + \b, \offsety+\x * \b + \b) rectangle (\offsetx+\m*\b - \x * \b , \offsety+\x * \b );
\foreach \x in {2,...,\m}
\draw[thick, white,fill=darkgray] (\offsetx +\m*\b - \x * \b + 2*\b, \offsety+\x * \b + \b) rectangle (\offsetx +\m*\b - \x * \b + \b, \offsety+\x * \b );
\foreach \x in {2,...,\m}
\draw[thick, white, fill=darkgray] (\offsetx +\m*\b - \x * \b + \b, \offsety+\x * \b) rectangle (\offsetx +\m*\b - \x * \b, \offsety+\x * \b - \b);

\foreach \x in {1,...,\m}
\draw[thick, white] (\offsetx + \x *\b,\offsety+\b) -- (\offsetx + \x * \b,\offsety+2*\m+\b);

\filldraw[fill=none, thick] (\offsetx +0,\offsety+\b) rectangle (\offsetx +\m*\b,\offsety+\m*\b+\b);

\pgfmathsetmacro{\offsetx}{7.5*\m}
\pgfmathsetmacro{\offsety}{0}

\draw[fill=darkgray] (\offsetx,\m+1.5*\b) rectangle (\offsetx+\b,\m+2.5*\b);
\draw[fill=gray] (\offsetx,\m ) rectangle (\offsetx+\b, \m+ 1.0*\b);
\node at (\offsetx+4*\b,\m+2*\b) {full rank};
\node at (\offsetx+4*\b,\m+0.5*\b) {low rank};

\end{tikzpicture}

%% file: FMM_typicalFMM_step1.tex
\begin{tikzpicture}[scale=0.14]

\pgfmathsetmacro{\b}{2}
\pgfmathsetmacro{\k}{1}

\pgfmathsetmacro{\m}{8}

\pgfmathsetmacro{\offsetx}{2.5*\m}
\pgfmathsetmacro{\offsety}{0}

\filldraw[fill=gray, thick] (\b,\b) rectangle (\m*\b+\b,\m*\b+\b);

\foreach \x in {1,...,\m}
\foreach \y in {1,...,\m}
\draw[thick] (\x*\b,\y*\b) rectangle (\x*\b + \b, \y*\b + \b);
\node at (\m+\b,-\b) {$\underset{N \times N}{\A}$};

\node at (2*\m+3*\b,\m + \b) {$\approx$};

\pgfmathsetmacro{\offsetx}{2.75*\m + 2*\b}
\pgfmathsetmacro{\offsety}{0}
\filldraw[fill=none, thick] (\offsetx +0,\offsety+\b) rectangle (\offsetx +\m*\b,\offsety+\m*\b+\b);

\foreach \x in {1,...,\m}
\draw[thick,fill=gray] (\offsetx + \m*\b - \x * \b + \b, \offsety+\x * \b + \b) rectangle (\offsetx+\m*\b - \x * \b , \offsety+\x * \b );
\foreach \x in {2,...,\m}
\draw[thick,fill=gray] (\offsetx +\m*\b - \x * \b + 2*\b, \offsety+\x * \b + \b) rectangle (\offsetx +\m*\b - \x * \b + \b, \offsety+\x * \b );
\foreach \x in {2,...,\m}
\draw[thick, fill=gray] (\offsetx +\m*\b - \x * \b + \b, \offsety+\x * \b) rectangle (\offsetx +\m*\b - \x * \b, \offsety+\x * \b - \b);

\node at (\offsetx+\m,-\b) {$\underset{N \times N}{\A^{\rm (near)}}$};
\node at (5*\m+3.5*\b,\m +\b) {$+$};

\pgfmathsetmacro{\offsetx}{5.25*\m+5*\b}
\pgfmathsetmacro{\offsety}{0}
\filldraw[fill=none, thick] (\offsetx +0,\offsety+\b) rectangle (\offsetx +\m*\k,\offsety+\m*\b+\b);
\foreach \x in {1,...,\m}
\draw[thick,fill=gray] (\offsetx + \m*\k - \x * \k + \k, \offsety+\x * \b + \b) rectangle (\offsetx+\m*\k - \x*\k, \offsety+\x * \b );
\node at (\offsetx+1.5*\b,-\b) {$\underset{N \times pk}{\mtx E}$};

\pgfmathsetmacro{\offsetx}{5.25*\m + \m*\k+ 6*\b}
\pgfmathsetmacro{\offsety}{0}
\filldraw[fill=gray, thick] (\offsetx +0,\offsety+\b + \m*\k) rectangle (\offsetx +\m*\k,\offsety+\m*\b+\b);

\pgfmathsetmacro{\offsety}{\m*\k+\k}
\foreach \x in {1,...,\m}
\draw[thick] (\offsetx, \offsety + \x * \k) -- (\offsetx + \m*\k, \offsety + \x * \k);
\foreach \x in {1,...,\m}
\draw[thick] (\offsetx + \x * \k, \offsety+\k) -- (\offsetx + \x * \k, \offsety + \m * \k+\k);

\foreach \x in {1,...,\m}
\draw[thick,fill=white] (\offsetx + \m*\k - \x * \k + \k, \offsety+\x * \k + \k) rectangle (\offsetx+\m*\k - \x * \k , \offsety+\x * \k );
\foreach \x in {2,...,\m}
\draw[thick,fill=white] (\offsetx +\m*\k - \x * \k + 2*\k, \offsety+\x * \k + \k) rectangle (\offsetx +\m*\k - \x * \k+ \k, \offsety+\x * \k );
\foreach \x in {2,...,\m}
\draw[thick, fill=white] (\offsetx +\m*\k - \x * \k + \k, \offsety+\x * \k) rectangle (\offsetx +\m*\k - \x * \k, \offsety+\x * \k - \k);
\node at (\offsetx+0.5*\m*\k,-\b) {$\underset{pk \times pk}{\widetilde{\mtx A}}$};

\pgfmathsetmacro{\offsetx}{5.25*\m + 2*\m*\k+ 7*\b}
\pgfmathsetmacro{\offsety}{0}
\filldraw[fill=none, thick] (\offsetx +0,\offsety+\b + \m*\k) rectangle (\offsetx +\m*\b,\offsety+\m*\b+\b);
\pgfmathsetmacro{\offsety}{\m*\k+\k}
\foreach \x in {1,...,\m}
\draw[thick,fill=gray] (\offsetx + \m*\b - \x * \b + \b, \offsety+\x * \k + \k) rectangle (\offsetx+\m*\b - \x*\b, \offsety+\x * \k );
\node at (\offsetx + 3*\b,-\b) {$\underset{pk \times N}{\mtx F}$};

\end{tikzpicture}

%% file: FMM_typicalFMM_step2prime.tex
\begin{tikzpicture}[scale=0.2]

\pgfmathsetmacro{\b}{2}
\pgfmathsetmacro{\k}{1}

\pgfmathsetmacro{\m}{8}

\pgfmathsetmacro{\offsetx}{\m*\k+ \k}
\pgfmathsetmacro{\offsety}{0}
\filldraw[fill=gray, thick] (\offsetx +0,\offsety) rectangle (\offsetx +\m*\k,\offsety+\m*\k);

\foreach \x in {1,...,\m}
\draw[thick] (\offsetx, \offsety + \x * \k) -- (\offsetx + \m*\k, \offsety + \x * \k);
\foreach \x in {1,...,\m}
\draw[thick] (\offsetx + \x * \k, \offsety) -- (\offsetx + \x * \k, \offsety + \m * \k);

\foreach \x in {1,...,\m}
\draw[thick,fill=white] (\offsetx + \m*\k - \x * \k + \k, \offsety+\x * \k-\k) rectangle (\offsetx+\m*\k - \x * \k , \offsety+\x * \k);
\foreach \x in {2,...,\m}
\draw[thick,fill=white] (\offsetx +\m*\k - \x * \k + 2*\k, \offsety+\x * \k) rectangle (\offsetx +\m*\k - \x * \k+ \k, \offsety+\x * \k -\k);
\foreach \x in {2,...,\m}
\draw[thick, fill=white] (\offsetx +\m*\k - \x * \k + \k, \offsety+\x * \k-2*\k) rectangle (\offsetx +\m*\k - \x * \k, \offsety+\x * \k - \k);

\node at (\offsetx+0.5*\m*\k,-1.5*\b) {$\underset{pk \times pk}{\widetilde{\mtx A}}$};
\node at (\offsetx+\m*\k+\b,+0.5*\m*\k) {$=$};

\pgfmathsetmacro{\offsetx}{3.5*\m*\k + \k + 2*\b}
\pgfmathsetmacro{\offsetx}{2.5*\m*\k+\k}
\pgfmathsetmacro{\offsety}{0}

\filldraw[fill=gray, thick] (\offsetx +0,\offsety) rectangle (\offsetx +\m*\k,\offsety+\m*\k);

\foreach \x in {1,...,\m}
\fill[fill=white] (\offsetx + \m*\k - \x * \k + \k, \offsety+\x * \k-\k) rectangle (\offsetx+\m*\k - \x * \k , \offsety+\x * \k);
\foreach \x in {2,...,\m}
\fill[fill=white] (\offsetx +\m*\k - \x * \k + 2*\k, \offsety+\x * \k) rectangle (\offsetx +\m*\k - \x * \k+ \k, \offsety+\x * \k -\k);
\foreach \x in {2,...,\m}
\fill[fill=white] (\offsetx +\m*\k - \x * \k + \k, \offsety+\x * \k-2*\k) rectangle (\offsetx +\m*\k - \x * \k, \offsety+\x * \k - \k); 

\pgfmathsetmacro{\iblock}{4}; \pgfmathsetmacro{\jblock}{3}
\fill[fill=white] (\offsetx + \iblock*2*\k - 2*\k, \offsety+\jblock *2*\k  - 2*\k) rectangle (\offsetx+\iblock*2*\k, \offsety+\jblock*2*\k);
\pgfmathsetmacro{\iblock}{4}; \pgfmathsetmacro{\jblock}{4}
\fill[fill=white] (\offsetx + \iblock*2*\k - 2*\k, \offsety+\jblock *2*\k  - 2*\k) rectangle (\offsetx+\iblock*2*\k, \offsety+\jblock*2*\k);
\pgfmathsetmacro{\iblock}{3}; \pgfmathsetmacro{\jblock}{4}
\fill[fill=white] (\offsetx + \iblock*2*\k - 2*\k, \offsety+\jblock *2*\k  - 2*\k) rectangle (\offsetx+\iblock*2*\k, \offsety+\jblock*2*\k);
\pgfmathsetmacro{\iblock}{1}; \pgfmathsetmacro{\jblock}{2}
\fill[fill=white] (\offsetx + \iblock*2*\k - 2*\k, \offsety+\jblock *2*\k  - 2*\k) rectangle (\offsetx+\iblock*2*\k, \offsety+\jblock*2*\k);
\pgfmathsetmacro{\iblock}{2}; \pgfmathsetmacro{\jblock}{1}
\fill[fill=white] (\offsetx + \iblock*2*\k - 2*\k, \offsety+\jblock *2*\k  - 2*\k) rectangle (\offsetx+\iblock*2*\k, \offsety+\jblock*2*\k);
\pgfmathsetmacro{\iblock}{1}; \pgfmathsetmacro{\jblock}{1}
\fill[fill=white] (\offsetx + \iblock*2*\k - 2*\k, \offsety+\jblock *2*\k  - 2*\k) rectangle (\offsetx+\iblock*2*\k, \offsety+\jblock*2*\k);

\foreach \x in {1,...,\m}
\draw[thick] (\offsetx, \offsety + \x * \k) -- (\offsetx + \m*\k, \offsety + \x * \k);
\foreach \x in {1,...,\m}
\draw[thick] (\offsetx + \x * \k, \offsety) -- (\offsetx + \x * \k, \offsety + \m * \k);

\node at (\offsetx+0.5*\m*\k,-\b) {${\mtx A}_{\S,\S}^{\rm (interaction\ list)}$};
\node at (\offsetx+\m*\k+\b,+0.5*\m*\k) {$+$};

\pgfmathsetmacro{\offsetx}{3.5*\m*\k + \k + 2*\b}
\pgfmathsetmacro{\offsety}{0}

\filldraw[fill=gray, thick] (\offsetx +0,\offsety) rectangle (\offsetx +\m*\k,\offsety+\m*\k);

\pgfmathsetmacro{\m}{\m*0.5}
\pgfmathsetmacro{\k}{\k*2}

\foreach \x in {1,...,\m}
\draw[thick] (\offsetx, \offsety + \x * \k) -- (\offsetx + \m*\k, \offsety + \x * \k);
\foreach \x in {1,...,\m}
\draw[thick] (\offsetx + \x * \k, \offsety) -- (\offsetx + \x * \k, \offsety + \m * \k);

\foreach \x in {1,...,\m}
\draw[thick,fill=white] (\offsetx + \m*\k - \x * \k + \k, \offsety+\x * \k-\k) rectangle (\offsetx+\m*\k - \x * \k , \offsety+\x * \k);
\foreach \x in {2,...,\m}
\draw[thick,fill=white] (\offsetx +\m*\k - \x * \k + 2*\k, \offsety+\x * \k) rectangle (\offsetx +\m*\k - \x * \k+ \k, \offsety+\x * \k -\k);
\foreach \x in {2,...,\m}
\draw[thick, fill=white] (\offsetx +\m*\k - \x * \k + \k, \offsety+\x * \k-2*\k) rectangle (\offsetx +\m*\k - \x * \k, \offsety+\x * \k - \k);

\node at (\offsetx+0.5*\m*\k,-\b) {${\mtx A}_{\S,\S}^{\rm (far)}$};

\end{tikzpicture}

%% file: FMM_modifiedFMM_step1.tex
\begin{tikzpicture}[scale=0.12]

\pgfmathsetmacro{\b}{2}
\pgfmathsetmacro{\k}{1}

\pgfmathsetmacro{\m}{8}

\pgfmathsetmacro{\offsetx}{2.5*\m}
\pgfmathsetmacro{\offsety}{0}

\filldraw[fill=gray, thick] (\b,\b) rectangle (\m*\b+\b,\m*\b+\b);

\foreach \x in {1,...,\m}
\foreach \y in {1,...,\m}
\draw[thick] (\x*\b,\y*\b) rectangle (\x*\b + \b, \y*\b + \b);
\node at (\m+\b,-1*\b) {$\underset{N \times N}{\A}$};

\node at (2*\m+3*\b,\m + \b) {$\approx$};

\pgfmathsetmacro{\offsetx}{2.75*\m + 2*\b}
\pgfmathsetmacro{\offsety}{0}
\filldraw[fill=none, thick] (\offsetx +0,\offsety+\b) rectangle (\offsetx +\m*\b,\offsety+\m*\b+\b);

\foreach \x in {1,...,\m}
\draw[thick,fill=gray] (\offsetx + \m*\b - \x * \b + \b, \offsety+\x * \b + \b) rectangle (\offsetx+\m*\b - \x * \b , \offsety+\x * \b );
\foreach \x in {2,...,\m}
\draw[thick,fill=gray] (\offsetx +\m*\b - \x * \b + 2*\b, \offsety+\x * \b + \b) rectangle (\offsetx +\m*\b - \x * \b + \b, \offsety+\x * \b );
\foreach \x in {2,...,\m}
\draw[thick, fill=gray] (\offsetx +\m*\b - \x * \b + \b, \offsety+\x * \b) rectangle (\offsetx +\m*\b - \x * \b, \offsety+\x * \b - \b);

\node at (\offsetx+\m,-1*\b) {$\underset{N \times N}{\A^{\rm (near)}}$};
\node at (5*\m+3.5*\b,\m +\b) {$-$};

\pgfmathsetmacro{\offsetx}{5.25*\m+5*\b}
\pgfmathsetmacro{\offsety}{0}
\filldraw[fill=none, thick] (\offsetx +0,\offsety+\b) rectangle (\offsetx +\m*\k,\offsety+\m*\b+\b);
\foreach \x in {1,...,\m}
\draw[thick,fill=gray] (\offsetx + \m*\k - \x * \k + \k, \offsety+\x * \b + \b) rectangle (\offsetx+\m*\k - \x*\k, \offsety+\x * \b );
\node at (\offsetx+1.5*\b,-1*\b) {$\underset{N \times pk}{\mtx E}$};

\pgfmathsetmacro{\offsetx}{\offsetx + \m*\k+ \b}
\pgfmathsetmacro{\offsety}{0}
\filldraw[fill=white, thick] (\offsetx +0,\offsety+\b + \m*\k) rectangle (\offsetx +\m*\k,\offsety+\m*\b+\b);

\pgfmathsetmacro{\offsety}{\m*\k+\k}
\foreach \x in {1,...,\m}
\draw[thick] (\offsetx, \offsety + \x * \k) -- (\offsetx + \m*\k, \offsety + \x * \k);
\foreach \x in {1,...,\m}
\draw[thick] (\offsetx + \x * \k, \offsety+\k) -- (\offsetx + \x * \k, \offsety + \m * \k+\k);

\foreach \x in {1,...,\m}
\draw[thick,fill=gray] (\offsetx + \m*\k - \x * \k + \k, \offsety+\x * \k + \k) rectangle (\offsetx+\m*\k - \x * \k , \offsety+\x * \k );
\foreach \x in {2,...,\m}
\draw[thick,fill=gray] (\offsetx +\m*\k - \x * \k + 2*\k, \offsety+\x * \k + \k) rectangle (\offsetx +\m*\k - \x * \k+ \k, \offsety+\x * \k );
\foreach \x in {2,...,\m}
\draw[thick, fill=gray] (\offsetx +\m*\k - \x * \k + \k, \offsety+\x * \k) rectangle (\offsetx +\m*\k - \x * \k, \offsety+\x * \k - \k);
\node at (\offsetx+0.5*\m*\k,-\b) {$\underset{pk \times pk} { {\mtx A}_{\S,\S}^{\rm (near)} }$};

\pgfmathsetmacro{\offsetx}{\offsetx+\m*\k+\b}
\pgfmathsetmacro{\offsety}{0}
\filldraw[fill=none, thick] (\offsetx +0,\offsety+\b + \m*\k) rectangle (\offsetx +\m*\b,\offsety+\m*\b+\b);
\pgfmathsetmacro{\offsety}{\m*\k+\k}
\foreach \x in {1,...,\m}
\draw[thick,fill=gray] (\offsetx + \m*\b - \x * \b + \b, \offsety+\x * \k + \k) rectangle (\offsetx+\m*\b - \x*\b, \offsety+\x * \k );
\node at (\offsetx + 3*\b,-\b) {$\underset{pk \times N}{\mtx F}$ };

\node at (10*\m + 6*\b,\m +\b) {$+$};

\pgfmathsetmacro{\offsetx}{10*\m+8*\b}
\pgfmathsetmacro{\offsety}{0}
\filldraw[fill=none, thick] (\offsetx +0,\offsety+\b) rectangle (\offsetx +\m*\k,\offsety+\m*\b+\b);
\foreach \x in {1,...,\m}
\draw[thick,fill=gray] (\offsetx + \m*\k - \x * \k + \k, \offsety+\x * \b + \b) rectangle (\offsetx+\m*\k - \x*\k, \offsety+\x * \b );
\node at (\offsetx+1.5*\b,-\b) {$\underset{N \times pk}{\mtx E}$};

\pgfmathsetmacro{\offsetx}{\offsetx + \m*\k+ \b}
\pgfmathsetmacro{\offsety}{0}
\filldraw[fill=gray, thick] (\offsetx +0,\offsety+\b + \m*\k) rectangle (\offsetx +\m*\k,\offsety+\m*\b+\b);

\pgfmathsetmacro{\offsety}{\m*\k+\k}
\foreach \x in {1,...,\m}
\draw[thick] (\offsetx, \offsety + \x * \k) -- (\offsetx + \m*\k, \offsety + \x * \k);
\foreach \x in {1,...,\m}
\draw[thick] (\offsetx + \x * \k, \offsety+\k) -- (\offsetx + \x * \k, \offsety + \m * \k+\k);

\node at (\offsetx+0.5*\m*\k,-\b) {$\underset{pk \times pk}{\mtx A_{\S,\S}}$};

\pgfmathsetmacro{\offsetx}{\offsetx+\m*\k+\b}
\pgfmathsetmacro{\offsety}{0}
\filldraw[fill=none, thick] (\offsetx +0,\offsety+\b + \m*\k) rectangle (\offsetx +\m*\b,\offsety+\m*\b+\b);
\pgfmathsetmacro{\offsety}{\m*\k+\k}
\foreach \x in {1,...,\m}
\draw[thick,fill=gray] (\offsetx + \m*\b - \x * \b + \b, \offsety+\x * \k + \k) rectangle (\offsetx+\m*\b - \x*\b, \offsety+\x * \k );
\node at (\offsetx + 3*\b,-\b) {$\underset{pk \times N}{\mtx F}$};

\end{tikzpicture}

%% file: FMM_modifiedFMM_step2.tex
\begin{tikzpicture}[scale=0.18]

\pgfmathsetmacro{\b}{2}
\pgfmathsetmacro{\k}{1}
\pgfmathsetmacro{\m}{8}

\pgfmathsetmacro{\offsetx}{\m*\k+ \k}
\pgfmathsetmacro{\offsety}{0}
\filldraw[fill=gray, thick] (\offsetx +0,\offsety) rectangle (\offsetx +\m*\k,\offsety+\m*\k);

\foreach \x in {1,...,\m}
\draw[thick] (\offsetx, \offsety + \x * \k) -- (\offsetx + \m*\k, \offsety + \x * \k);
\foreach \x in {1,...,\m}
\draw[thick] (\offsetx + \x * \k, \offsety) -- (\offsetx + \x * \k, \offsety + \m * \k);

\node at (\offsetx+0.5*\m*\k,-1.5*\b) {$\underset{pk \times pk}{{\mtx A}_{\S,\S}}$};
\node at (\offsetx+\m*\k+\b,+0.5*\m*\k) {$\approx$};

\pgfmathsetmacro{\b}{2}
\pgfmathsetmacro{\k}{1}
\pgfmathsetmacro{\m}{4}

\pgfmathsetmacro{\offsetx}{4.5*\m + 2*\b}
\pgfmathsetmacro{\offsety}{-0.5*\m}
\filldraw[fill=none, thick] (\offsetx +0,\offsety+\b) rectangle (\offsetx +\m*\b,\offsety+\m*\b+\b);

\foreach \x in {1,...,\m}
\draw[thick,fill=gray] (\offsetx + \m*\b - \x * \b + \b, \offsety+\x * \b + \b) rectangle (\offsetx+\m*\b - \x * \b , \offsety+\x * \b );
\foreach \x in {2,...,\m}
\draw[thick,fill=gray] (\offsetx +\m*\b - \x * \b + 2*\b, \offsety+\x * \b + \b) rectangle (\offsetx +\m*\b - \x * \b + \b, \offsety+\x * \b );
\foreach \x in {2,...,\m}
\draw[thick, fill=gray] (\offsetx +\m*\b - \x * \b + \b, \offsety+\x * \b) rectangle (\offsetx +\m*\b - \x * \b, \offsety+\x * \b - \b);

\node at (8*\m,0.5*\m +\b) {$-$};

\pgfmathsetmacro{\offsetx}{6*\m+5*\b}
\pgfmathsetmacro{\offsety}{-0.5*\m}
\filldraw[fill=none, thick] (\offsetx +0,\offsety+\b) rectangle (\offsetx +\m*\k,\offsety+\m*\b+\b);
\foreach \x in {1,...,\m}
\draw[thick,fill=gray] (\offsetx + \m*\k - \x * \k + \k, \offsety+\x * \b + \b) rectangle (\offsetx+\m*\k - \x*\k, \offsety+\x * \b );

\pgfmathsetmacro{\offsetx}{\offsetx + \m*\k+ \k}
\pgfmathsetmacro{\offsety}{-0.5*\m}
\filldraw[fill=white, thick] (\offsetx +0,\offsety+\b + \m*\k) rectangle (\offsetx +\m*\k,\offsety+\m*\b+\b);

\pgfmathsetmacro{\offsety}{-0.5*\m + \m*\k+\k}
\foreach \x in {1,...,\m}
\draw[thick] (\offsetx, \offsety + \x * \k) -- (\offsetx + \m*\k, \offsety + \x * \k);
\foreach \x in {1,...,\m}
\draw[thick] (\offsetx + \x * \k, \offsety+\k) -- (\offsetx + \x * \k, \offsety + \m * \k+\k);

\foreach \x in {1,...,\m}
\draw[thick,fill=gray] (\offsetx + \m*\k - \x * \k + \k, \offsety+\x * \k + \k) rectangle (\offsetx+\m*\k - \x * \k , \offsety+\x * \k );
\foreach \x in {2,...,\m}
\draw[thick,fill=gray] (\offsetx +\m*\k - \x * \k + 2*\k, \offsety+\x * \k + \k) rectangle (\offsetx +\m*\k - \x * \k+ \k, \offsety+\x * \k );
\foreach \x in {2,...,\m}
\draw[thick, fill=gray] (\offsetx +\m*\k - \x * \k + \k, \offsety+\x * \k) rectangle (\offsetx +\m*\k - \x * \k, \offsety+\x * \k - \k);

\pgfmathsetmacro{\offsetx}{\offsetx+\m*\k+\k}
\pgfmathsetmacro{\offsety}{-0.5*\m}
\filldraw[fill=none, thick] (\offsetx +0,\offsety+\b + \m*\k) rectangle (\offsetx +\m*\b,\offsety+\m*\b+\b);
\pgfmathsetmacro{\offsety}{-0.5*\m + \m*\k+\k}
\foreach \x in {1,...,\m}
\draw[thick,fill=gray] (\offsetx + \m*\b - \x * \b + \b, \offsety+\x * \k + \k) rectangle (\offsetx+\m*\b - \x*\b, \offsety+\x * \k );

\node at (10*\m + 7*\b,0.5*\m +\b) {$+$};

\pgfmathsetmacro{\offsetx}{10*\m+8*\b}
\pgfmathsetmacro{\offsety}{-0.5*\m}
\filldraw[fill=none, thick] (\offsetx +0,\offsety+\b) rectangle (\offsetx +\m*\k,\offsety+\m*\b+\b);
\foreach \x in {1,...,\m}
\draw[thick,fill=gray] (\offsetx + \m*\k - \x * \k + \k, \offsety+\x * \b + \b) rectangle (\offsetx+\m*\k - \x*\k, \offsety+\x * \b );

\pgfmathsetmacro{\offsetx}{\offsetx + \m*\k+ \k}
\pgfmathsetmacro{\offsety}{-0.5*\m}
\filldraw[fill=gray, thick] (\offsetx +0,\offsety+\b + \m*\k) rectangle (\offsetx +\m*\k,\offsety+\m*\b+\b);

\pgfmathsetmacro{\offsety}{-0.5*\m + \m*\k+\k}
\foreach \x in {1,...,\m}
\draw[thick] (\offsetx, \offsety + \x * \k) -- (\offsetx + \m*\k, \offsety + \x * \k);
\foreach \x in {1,...,\m}
\draw[thick] (\offsetx + \x * \k, \offsety+\k) -- (\offsetx + \x * \k, \offsety + \m * \k+\k);


\pgfmathsetmacro{\offsetx}{\offsetx+\m*\k+\k}
\pgfmathsetmacro{\offsety}{-0.5*\m}
\filldraw[fill=none, thick] (\offsetx +0,\offsety+\b + \m*\k) rectangle (\offsetx +\m*\b,\offsety+\m*\b+\b);
\pgfmathsetmacro{\offsety}{-0.5*\m + \m*\k+\k}
\foreach \x in {1,...,\m}
\draw[thick,fill=gray] (\offsetx + \m*\b - \x * \b + \b, \offsety+\x * \k + \k) rectangle (\offsetx+\m*\b - \x*\b, \offsety+\x * \k );

\end{tikzpicture}

%% file: alg_fact.tex
\begin{algorithm}[!htb]
\caption{\texttt{Strong recursive skeletonization}}
\label{a:precompute}
\begin{algorithmic}[1]
\Require{Kernel function $K$ and tree decomposition $\T$ of points $\X$.}
\Ensure A partition $\B_i = \S_i \cup \R_i$ and an interpolation matrix $\mtx T_{\B_i}$ for each box $\B_i$.
\Statex
\For{$\ell = L, L-1, \ldots, 1$} 
\ForAll {boxes $\B_i$ at level $\ell$}
\Comment{\hfill \parbox[t]{.4\linewidth}{$\triangleright$ In parallel at each level.}}
\If {box $\B_i$ is a leaf}
\State Determine the index vector $\B_i$ from the tree decomposition $\T$.
\Else 
\State Form the index vector $\B_i$ by aggregating the skeleton indices of the children as \[\B_i := \bigcup_{\B_j \in \C(\B_i)} {\S_j}.\]
\EndIf
\State{Form a partition $\B_i  = \S_i \cup \R_i$ and a matrix $\mtx T_{\B_i}$ using column ID of the proxy matrix.}
\EndFor
\EndFor

\end{algorithmic}
\end{algorithm}

%% file: alg_solve.tex
\clearpage
\begin{algorithm}[p]
\caption{\texttt{FMM apply}}
\label{a:evaluate}

\vspace*{\fill} 

\begin{algorithmic}[1]
\Require{$\vct{q} \in \mathbb{C}^N$. 
For each box $\B_i$, indices $I_{\B_i} = I_{\S_i} \cup I_{\R_i}$,  and matrix $\mtx{T}_{\B_i}$.}
\Ensure $\vct{u} = \A \, \vct{q}$.

\vspace{0.5em}

\ForAll{leaf boxes $\L_i$}
\State Set ${\mtx q}_{\L_i}$ to given charges.
\EndFor
\ForAll{boxes $\B_i$}
\State Set ${\mtx u}_{\B_i} := \mtx 0$ and $\hat{\mtx u}_{\B_i} := \mtx 0$.
\EndFor

\vspace{0.5em}

\For{$\ell = L, L-1, \ldots, 1$}  
\Comment{\hfill \parbox[t]{.3\linewidth}{$\triangleright$ Upward pass.}}
\ForAll{box $\B$ at level $\ell$} \Comment{\hfill \parbox[t]{.3\linewidth}{$\triangleright$ In parallel at each level.}}

\If{box $\B_i$ is a tree box}
\State Accumulate $\mtx q_{\B_i}$ from children as $$\mtx q_{\B_i}:= \left [ \hat{\mtx q}_{\S_j}; \hat{\mtx q}_{\S_{j+1}}; \dots \right]\ \text{for}\ \B_j, \B_{j+1}, \dots \in \C(\B_i).$$
\EndIf

\State Set outgoing representation
\[
\hat{\mtx q}_{\S_i}:= \begin{pmatrix} \mtx I & \mtx T_{\B_i}\end{pmatrix}
\begin{pmatrix}\mtx q_{\S_i}\\ \mtx q_{\R_i}\end{pmatrix}
\]
\EndFor
\EndFor


\vspace{0.5em}

\ForAll{boxes $\B_i$ with colleague neighbor $\B_j$} 
\Comment{\hfill \parbox[t]{.3\linewidth}{$\triangleright$ In parallel for all boxes.}}
\State Evaluate \Comment{\hfill \parbox[t]{.3\linewidth}{$\hphantom{\triangleright}$ Translation $\B_i \xleftarrow{\rm (ifo)} \B_j$.}}
\[
{\mtx u}_{\B_i} \pluseq \mtx A_{\B_i, \B_j} {\mtx q}_{\B_j}, \qquad
 \hat {\mtx u}_{\S_i} \minuseq \mtx A_{\S_i, \S_j} \hat {\mtx q}_{\S_j}
\]
\EndFor

\ForAll{boxes $\B_i$ with coarse neighbor $\L_j$} 
\Comment{\hfill \parbox[t]{.3\linewidth}{$\triangleright$ In parallel for all boxes.}}
\State Evaluate \Comment{\hfill \parbox[t]{.3\linewidth}{$\hphantom{\triangleright}$ Translation $\B_i \xleftarrow{\rm (ifs)} \L_j$.}}
\[
{\mtx u}_{\B_i} \pluseq \mtx A_{\B_i, \L_j} {\mtx q}_{\L_j}, \qquad
 \hat {\mtx u}_{\S_i} \minuseq \mtx A_{\S_i, \L_j} {\mtx q}_{\L_j}
\]
\EndFor

\ForAll{leaves $\L_i$ with fine neighbor $\B_j$} 
\Comment{\hfill \parbox[t]{.3\linewidth}{$\triangleright$ In parallel for all leaves.}}
\State Evaluate \Comment{\hfill \parbox[t]{.3\linewidth}{$\hphantom{\triangleright}$ Translation $\L_i \xleftarrow{\rm (tfo)} \B_j$.}}
\[
{\mtx u}_{\L_i} \pluseq \mtx A_{\L_i, \B_j} {\mtx q}_{\B_j} - \mtx A_{\L_i, \S_j} \hat {\mtx q}_{\S_j} 
\]
\EndFor

\ForAll{boxes $\B_i$ on level 1}
\State Evaluate \[{\mtx u}_{\B_i} \pluseq \sum_{\B_j\ \text{on level 1}} \mtx A_{\B_i, \B_j} {\mtx q}_{\B_j}.\]
\EndFor

\vspace{0.5em}



\For{$\ell = 2, \ldots, L$}  \Comment{\hfill \parbox[t]{.3\linewidth}{$\triangleright$ Downward pass.}}
\ForAll{boxes $\B_i$ at level $\ell$} \Comment{\hfill \parbox[t]{.3\linewidth}{$\triangleright$ In parallel at each level.}}
\State Add to $\hat{\mtx u}_{\S_i}$ from relevant indices of parent box $\mtx u_{\P(\B_i)}$ as $$\hat {\mtx{u}}_{\mathcal S_i} \pluseq \mtx u_{\mathcal S_p},\ \text{where}\ p = \mathcal P(\mathcal B_i).$$
\State Add to incoming representation $$\mtx u_{\B_i} \pluseq \begin{pmatrix} \mtx I\\ \mtx T_{\B_i}^* \end{pmatrix} \hat{\mtx u}_{\S_i} .$$
\EndFor
\EndFor

\ForAll{leaf boxes $\L_i$}
\State Set ${\mtx u}_{\L_i}$ to result vector $\mtx u$.
\EndFor
\end{algorithmic}
\vspace*{\fill} 
\end{algorithm}

\clearpage

%% file: results_2d_lp.tex
\clearpage

\noindent 
\vspace*{\fill} 

\begin{figure}[htbp]

\centering

\begin{subfigure}{0.5\textwidth}
\centering
\includegraphics[width=0.8\textwidth]{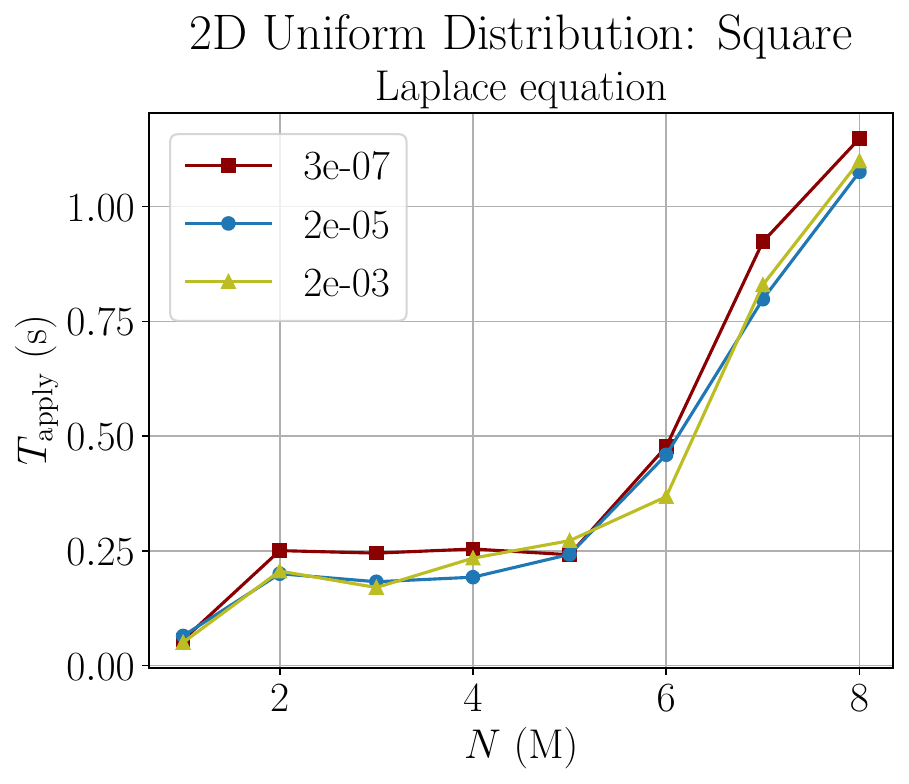}
\caption{2D Laplace kernel for points in square.}
\label{f:square_lp}
\end{subfigure}%
\begin{subfigure}{0.5\textwidth}
\centering
\includegraphics[width=0.85\textwidth]{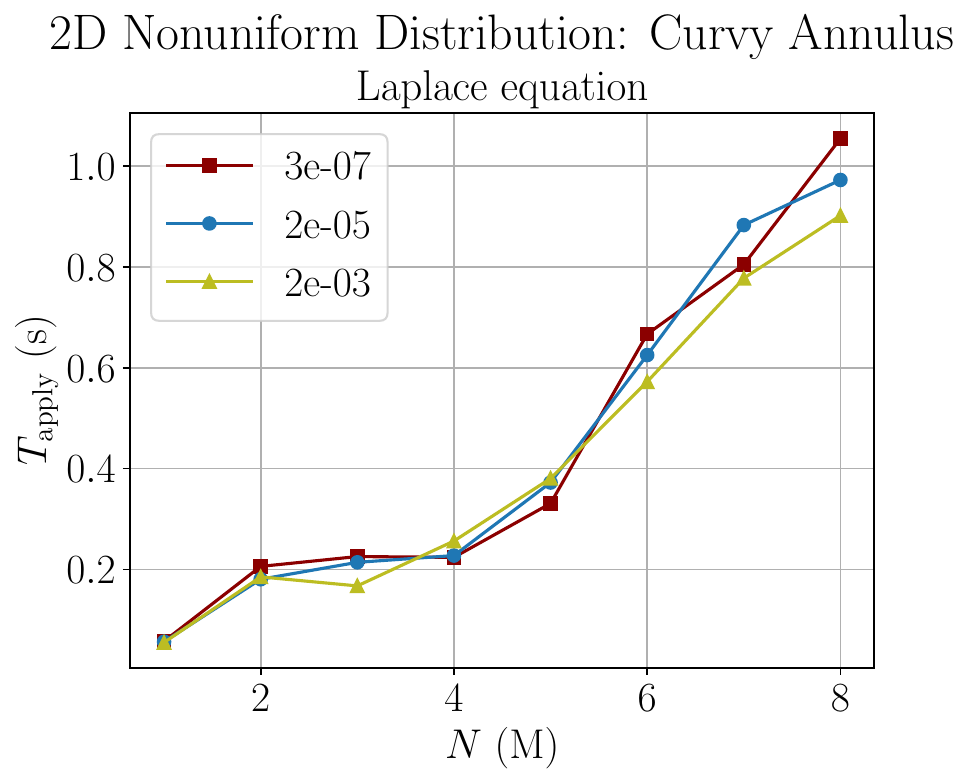}
\caption{2D Laplace kernel for points in curvy annulus.}
\label{f:curvy_lp}
\end{subfigure}

\vspace{0.5em}

\begin{subfigure}{\textwidth}
\centering \small
\begin{tabular}{c|cc|ccc|cc}
$N$ & $b$ & $k_{\rm max}$ & $T_{\rm tree}$ & $T_{\rm skel}$ & $M_{\rm proj}$ & $T_{\rm apply}$ & relerr \\ \hline
    &     & 10 &       & 4.0 s & 0.1 GB & 50 ms & 1.31e-03 \\
1 M & 100 & 18 & 1.8 s & 3.0 s & 0.2 GB & 65 ms & 7.32e-06 \\
    &     & 28 &       & 3.4 s & 0.2 GB & 56 ms & 2.26e-07 \\
 \hline
    &     & 10 &        & 69.2 s & 1.5 GB & 1099 ms & 9.36e-04 \\
8 M & 100 & 18 & 29.4 s & 56.3 s & 2.6 GB & 1076 ms & 8.70e-06 \\
    &     & 28 &        & 56.5 s & 3.8 GB & 1148 ms & 3.69e-07 \\
\end{tabular}
\caption{2D Laplace kernel for points in square.}
\label{t:square_lp}
\end{subfigure}

\vspace{0.5em}

\begin{subfigure}{\textwidth}
\centering \small
\begin{tabular}{c|cc|ccc|cc}
$N$ & $b$ & $k_{\rm max}$ & $T_{\rm tree}$ & $T_{\rm skel}$ & $M_{\rm proj}$ & $T_{\rm apply}$ & relerr \\ \hline
    &     & 9 &        & 3.5 s & 0.1 GB & 55 ms & 1.12e-03 \\
1 M & 100 & 17 & 4.5 s & 3.4 s & 0.2 GB & 57 ms & 1.96e-05 \\
    &     & 26 &       & 3.4 s & 0.3 GB & 58 ms & 4.45e-07 \\
 \hline
    &     & 9 &         & 55.6 s & 1.4 GB & 902 ms & 1.54e-03 \\
8 M & 100 & 17 & 52.2 s & 54.3 s & 2.5 GB & 972 ms & 1.01e-05 \\
    &     & 27 &        & 41.2 s & 3.6 GB & 1055 ms & 3.89e-07 \\
\end{tabular}
\caption{2D Laplace kernel for points in curvy annulus.}
\label{t:curvy_lp}
\end{subfigure}

\caption{The performance for the Laplace 2D kernel (\ref{e:laplace}) is demonstrated for two point distributions.
Figure \ref{f:square_lp} reports $T_{\rm apply}$ for a random distribution of points in $[0,1]^2$ with additional data reported in Table \ref{t:square_lp}.
Similarly, Figure \ref{f:curvy_lp} and Table \ref{t:curvy_lp} report results for points in a wiggly torus, as shown in Figure \ref{f:wiggly_geom}.
}
\label{f:2dlp}

\end{figure}

\vspace*{\fill} 

\clearpage

%% file: results_2d_hh.tex
\clearpage

\noindent 
\vspace*{\fill} 

\begin{figure}[htbp]

\centering
\begin{subfigure}{0.5\textwidth}
\centering
\includegraphics[width=0.8\textwidth]{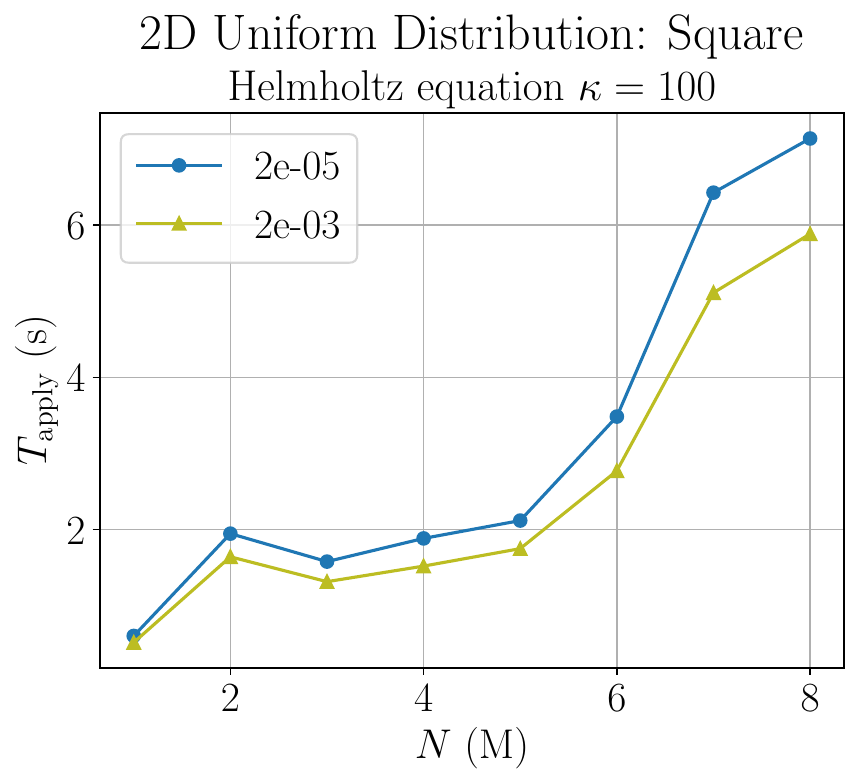}
\caption{2D Helmholtz kernel 
    for points in square.}
\label{f:square_hh}
\end{subfigure}%
\begin{subfigure}{0.5\textwidth}
\centering
\includegraphics[width=0.85\textwidth]{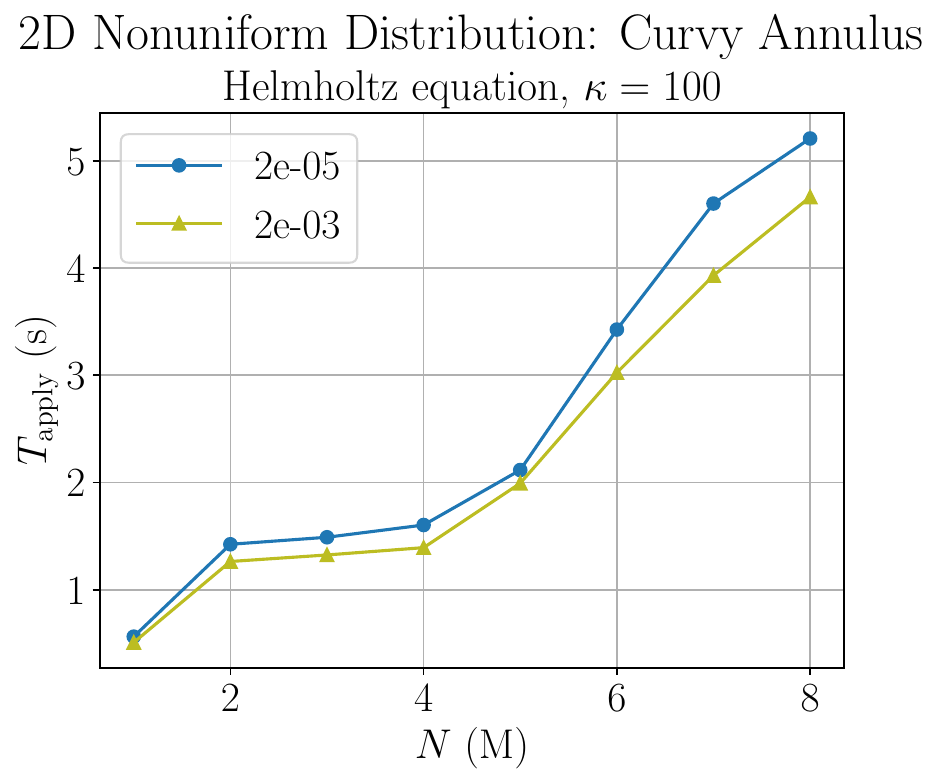}
\caption{ 2D Helmholtz kernel for points in curvy annulus.}
\label{f:curvy_hh}
\end{subfigure}

\vspace{0.5em}

\begin{subfigure}{\textwidth}
\centering \small
\begin{tabular}{c|cc|ccc|cc}
$N$ & $b$ & $k_{\rm max}$ & $T_{\rm tree}$ & $T_{\rm skel}$ & $M_{\rm proj}$ & $T_{\rm apply}$ & relerr \\ \hline
\multirow{2}{*}{1 M} & \multirow{2}{*}{100} & 32 & \multirow{2}{*}{2.0 s} & 3.3 s & 0.3 GB & 518 ms & 8.19e-04 \\
&  & 39 &  & 4.4 s & 0.5 GB & 604 ms & 7.73e-06 \\
 \hline
\multirow{2}{*}{8 M} & \multirow{2}{*}{100} & 32 & \multirow{2}{*}{30.1 s}  & 32.3 s & 5.1 GB & 5884 ms & 2.25e-03 \\
 &  & 40 & & 37.8 s & 8.3 GB & 7137 ms & 7.30e-06 \\
\end{tabular}
\caption{2D Helmholtz kernel ($\kappa=100$) for points in square.}
\label{t:square_hh}
\end{subfigure}

\vspace{0.5em}

\begin{subfigure}{\textwidth}
\centering \small
\begin{tabular}{c|cc|ccc|cc}
$N$ & $b$ & $k_{\rm max}$ & $T_{\rm tree}$ & $T_{\rm skel}$ & $M_{\rm proj}$ & $T_{\rm apply}$ & relerr \\ \hline
\multirow{2}{*}{1 M} & \multirow{2}{*}{100} & 29 & \multirow{2}{*}{4.4 s} & 3.1 s & 0.3 GB & 508 ms & 1.77e-03 \\
                     &  & 37 &  & 3.6 s & 0.5 GB & 562 ms & 1.14e-05 \\
\hline
\multirow{2}{*}{8 M} & \multirow{2}{*}{100} & 30 & \multirow{2}{*}{51.6 s} & 32.4 s & 3.8 GB & 4662 ms & 1.27e-03 \\
 &  & 36 & & 34.9 s & 6.6 GB & 5208 ms & 6.15e-06 \\
\end{tabular}
\caption{ 2D Helmholtz kernel ($\kappa=100$) for points in curvy annulus.}
\label{t:curvy_hh}
\end{subfigure}

\caption{The performance for the Helmholtz 2D kernel (\ref{e:helmholtz}) for $\kappa=100$ is demonstrated for two point distributions.
Figure \ref{f:square_hh} reports $T_{\rm apply}$ for a random distribution of points in $[0,1]^2$ with additional data reported in Table \ref{t:square_hh}.
Similarly, Figure \ref{f:curvy_hh} and Table \ref{t:curvy_hh} report results for points in a wiggly torus.
The use of single-precision with a complex-valued kernel leads to some loss of accuracy for the strictest compression tolerance,
we are limited to approximately 5 digits of accuracy, instead of 7, as for the 2D Laplace kernel, in Figure \ref{f:2dlp}.
}
\label{f:2dhh}
\end{figure}

\vspace*{\fill} 

\clearpage

%% file: results_3d_lp.tex
\clearpage

\noindent 
\vspace*{\fill} 

\begin{figure}[htb!]
\centering
\begin{subfigure}{0.5\textwidth}
\centering
\includegraphics[width=0.80\textwidth]{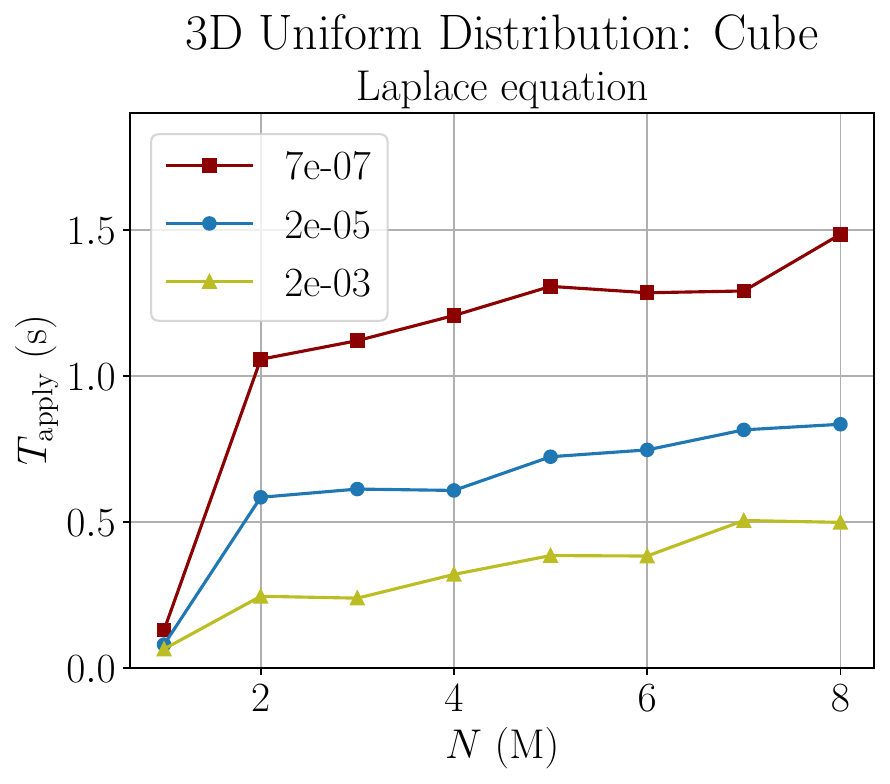}
\caption{3D Laplace kernel for points in cube.}
\label{f:cube_lp}
\end{subfigure}%
\begin{subfigure}{0.5\textwidth}
\centering
\includegraphics[width=0.85\textwidth]{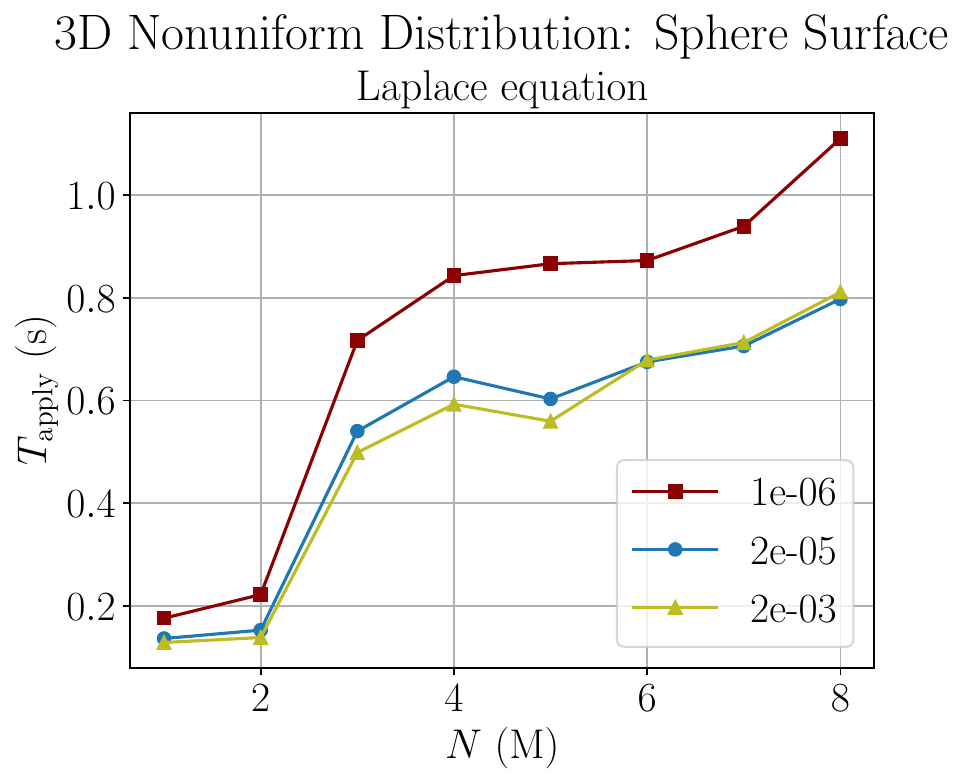}
\caption{3D Laplace kernel for points on sphere surface.}
\label{f:sphere_lp}
\end{subfigure}

\vspace{0.5em}

\begin{subfigure}{\textwidth}
\centering \small
\begin{tabular}{c|cc|ccc|cc}
$N$ & $b$ & $k_{\rm max}$ & $T_{\rm tree}$ & $T_{\rm skel}$ & $M_{\rm proj}$ & $T_{\rm apply}$ & relerr \\ \hline
    &     & 30 &       & 2.8 s & 0.2 GB & 64 ms & 1.34e-03 \\
1 M & 320 & 97 & 1.2 s & 6.1 s & 0.5 GB & 79 ms & 1.29e-05 \\
    &     & 181 &      & 14.9 s & 1.0 GB & 129 ms & 3.78e-07 \\
 \hline
    &     & 31 &         & 13.0 s & 1.4 GB & 498 ms & 1.76e-03 \\
8 M & 320 & 101 & 14.5 s & 27.6 s & 4.4 GB & 834 ms & 1.12e-05 \\
    &     & 190 &        & 60.7 s & 9.2 GB & 1485 ms & 7.51e-07 \\
\end{tabular}
\caption{3D Laplace kernel for points in cube.}
\label{t:cube_lp}
\end{subfigure}

\vspace{0.5em}

\begin{subfigure}{\textwidth}
\centering \small
\begin{tabular}{c|cc|ccc|cc}
$N$ & $b$ & $k_{\rm max}$ & $T_{\rm tree}$ & $T_{\rm skel}$ & $M_{\rm proj}$ & $T_{\rm apply}$ & relerr \\  \hline
    &     & 17 &       & 3.9 s & 0.3 GB & 129 ms & 2.64e-03 \\
1 M & 200 & 41 & 2.9 s & 5.0 s & 0.7 GB & 137 ms & 1.43e-05 \\
    &     & 75 &       & 7.2 s & 1.5 GB & 176 ms & 6.62e-07 \\
 \hline
    &     & 17 &        & 17.8 s & 1.7 GB & 811 ms & 1.32e-03 \\
8 M & 200 & 41 & 39.0 s & 22.4 s & 4.0 GB & 797 ms & 1.42e-05 \\
    &     & 75 &        & 30.2 s & 8.3 GB & 1110 ms & 1.26e-06 \\
\end{tabular}
\caption{3D Laplace kernel for points on sphere surface.}
\label{t:sphere_lp}
\end{subfigure}

\caption{The performance for the Laplace 3D kernel (\ref{e:laplace}) is demonstrated for two point distributions.
Figure \ref{f:cube_lp} reports $T_{\rm apply}$ for a random distribution of points in $[0,1]^3$ with additional data reported in Table \ref{t:cube_lp}.
Similarly, Figure \ref{f:sphere_lp} and Table \ref{t:sphere_lp} report results for points on the surface of the sphere.
}
\label{f:3dlp}
\end{figure}

\vspace*{\fill} 
\clearpage

%% file: results_3d_hh.tex
\clearpage

\noindent 
\vspace*{\fill} 

\begin{figure}[!htb]
\centering

\begin{subfigure}{0.5\textwidth}
\centering
\includegraphics[width=0.8\textwidth]{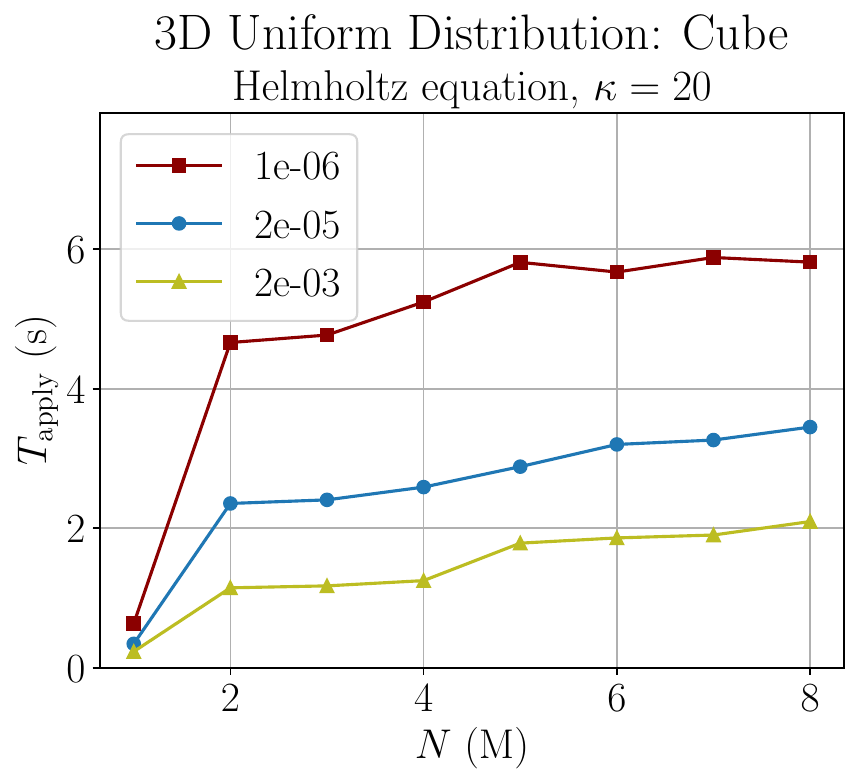}
\caption{3D Helmholtz kernel for points in cube.}
\label{f:cube_hh}
\end{subfigure}%
\begin{subfigure}{0.5\textwidth}
\centering
\includegraphics[width=0.85\textwidth]{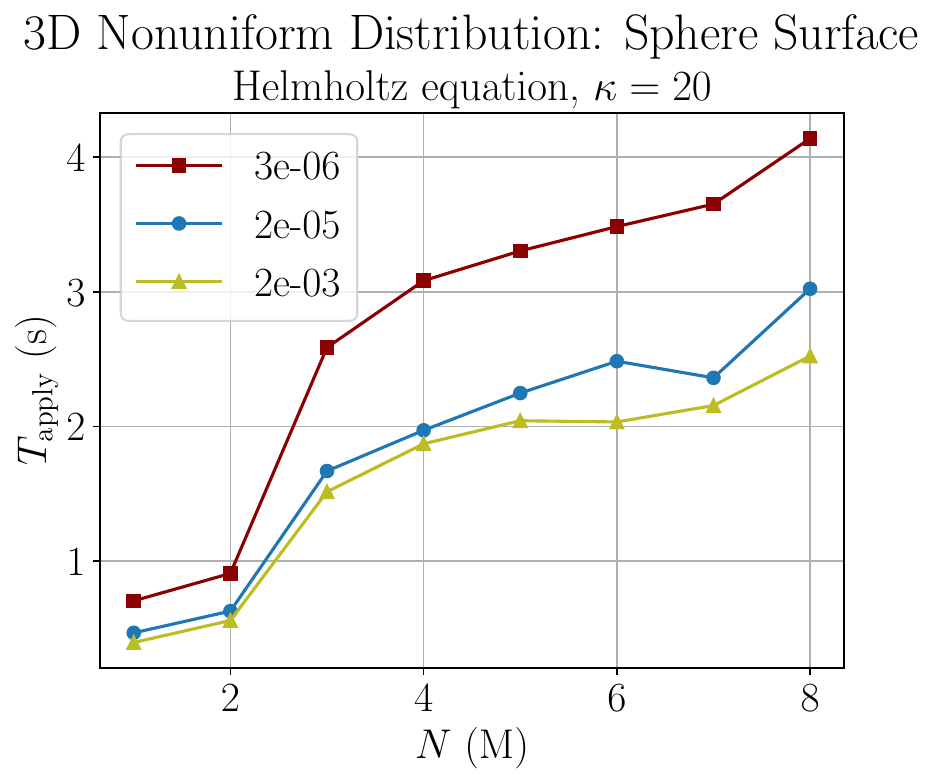}
\caption{3D Helmholtz kernel for points on sphere surface.}
\label{f:sphere_hh}
\end{subfigure}

\vspace{0.5em}

\begin{subfigure}{\textwidth}
\centering \small
\begin{tabular}{c|cc|ccc|cc}
$N$ & $b$ & $k_{\rm max}$ & $T_{\rm tree}$ & $T_{\rm skel}$ & $M_{\rm proj}$ & $T_{\rm apply}$ & relerr \\  \hline
    &     & 56 &        & 5.9 s & 0.5 GB & 231 ms & 9.56e-04 \\
1 M & 320 & 127 & 1.1 s & 13.3 s & 1.3 GB & 342 ms & 5.04e-06 \\
    &     & 211 &       & 28.8 s & 2.5 GB & 633 ms & 2.43e-06 \\
 \hline
    &     & 58 &         & 24.9 s & 3.4 GB & 2096 ms & 9.62e-04 \\
8 M & 320 & 133 & 14.6 s & 53.9 s & 10.9 GB & 3450 ms & 1.11e-05 \\
    &     & 225 &        & 133.8 s & 22.5 GB & 5816 ms & 1.60e-06 \\
\end{tabular}
\caption{3D Helmholtz kernel ($\kappa=20$) for points in cube.}
\label{t:cube_hh}
\end{subfigure}

\vspace{0.5em}

\begin{subfigure}{\textwidth}
\centering \small
\begin{tabular}{c|cc|ccc|cc}
$N$ & $b$ & $k_{\rm max}$ & $T_{\rm tree}$ & $T_{\rm skel}$ & $M_{\rm proj}$ & $T_{\rm apply}$ & relerr \\  \hline
    &     & 30 &       & 5.6 s & 0.7 GB & 397 ms & 1.04e-03 \\
1 M & 200 & 59 & 3.0 s & 7.7 s & 1.7 GB & 468 ms & 5.93e-06 \\
    &     & 96 &       & 12.1 s & 3.4 GB & 706 ms & 4.12e-06 \\
\hline
    &     & 29 &        & 27.0 s & 3.9 GB & 2524 ms & 1.50e-03 \\
8 M & 200 & 60 & 39.3 s & 35.0 s & 9.1 GB & 3023 ms & 1.05e-05 \\
    &     & 96 &        & 51.9 s & 19.0 GB & 4139 ms & 3.40e-06 \\
\end{tabular}
\caption{3D Helmholtz kernel ($\kappa=20$) for points on sphere surface.}
\label{t:sphere_hh}
\end{subfigure}
\caption{The performance for the Helmholtz 3D kernel (\ref{e:helmholtz}) for $\kappa=20$ is demonstrated for two point distributions.
Figure \ref{f:cube_hh} reports $T_{\rm apply}$ for a random distribution of points in $[0,1]^3$ with additional data reported in Table \ref{t:cube_hh}.
Similarly, Figure \ref{f:sphere_hh} and Table \ref{t:sphere_hh} report results for points on the surface of a sphere.}
\label{f:3dhh}
\end{figure}

\vspace*{\fill} 
\clearpage

%% file: results_comparison.tex
\clearpage

\noindent 
\vspace*{\fill} 

\begin{figure}[!htb]
\centering
\begin{subfigure}{0.5\textwidth}
\centering
\includegraphics[width=0.8\textwidth]{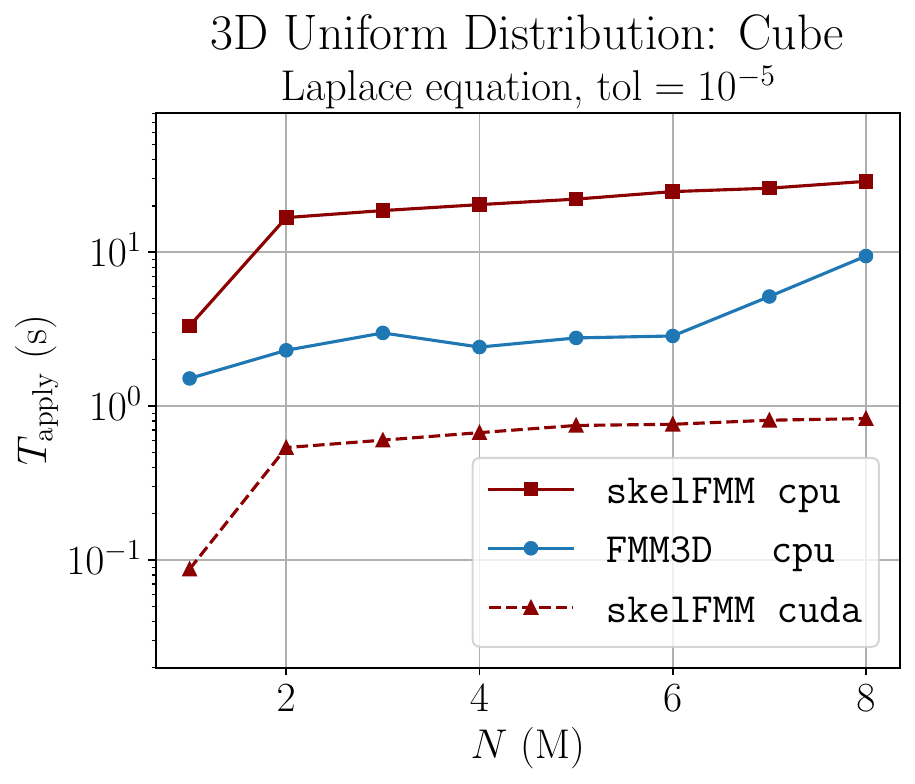}
\caption{Laplace kernel for points in a cube.}
\label{f:comp_lp_cube}
\end{subfigure}%
\begin{subfigure}{0.5\textwidth}
\centering
\includegraphics[width=0.85\textwidth]{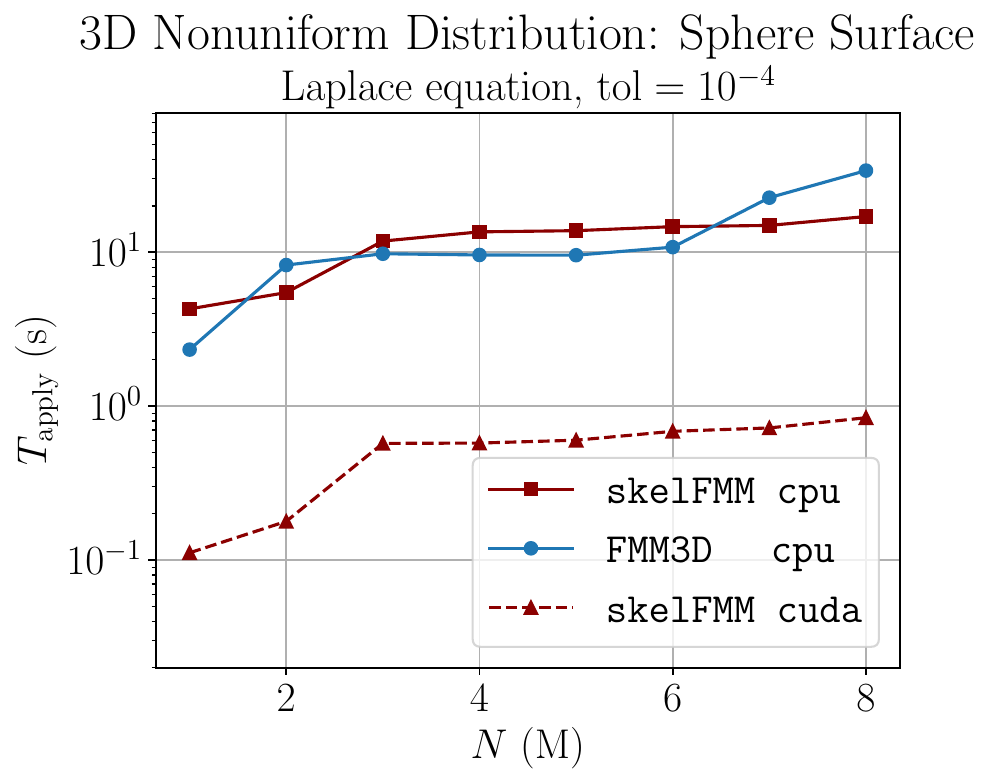}
\caption{Laplace kernel for points on sphere surface.}
\label{f:comp_lp_sphere}
\end{subfigure}

\vspace{0.5em}

\begin{subfigure}{0.5\textwidth}
\centering \small
\begin{tabular}{c|cc|cc}
& $N$ & $b$ & $T_{\rm apply}$ & relerr \\  \hline
\texttt{skelFMM cpu} &                            &  320  & 3.3 s & 1.28e-05 \\
\texttt{FMM3D}\hphantom{abc} \texttt{cpu}   & 1 M &  40   & 1.5 s & 2.65e-05 \\
\texttt{skelFMM cuda}&                            & 320   & 87 ms & 1.28e-05 \\
 \hline
\texttt{skelFMM cpu} &                            &  320  & 28.9 s  & 1.19e-05\\
\texttt{FMM3D}\hphantom{abc} \texttt{cpu}   & 8 M &  40   & 9.4 s  & 2.44e-05 \\
\texttt{skelFMM cuda}&                            & 320   & 831 ms & 1.19e-05\\
\end{tabular}
\caption{Laplace kernel for points in a cube.}
\label{t:comp_lp_cube}
\end{subfigure}%
\begin{subfigure}{0.5\textwidth}
\centering \small
\begin{tabular}{c|cc|cc}
& $N$ & $b$ & $T_{\rm apply}$ & relerr \\  \hline
\texttt{skelFMM cpu} &                            & 320   & 4.2 s & 9.67e-05\\
\texttt{FMM3D}\hphantom{abc} \texttt{cpu}   & 1 M &  40   & 2.3 s & 1.35e-04 \\
\texttt{skelFMM cuda}&                            &  320  & 112 ms & 9.56e-05\\
 \hline
\texttt{skelFMM cpu} &                            &  320  & 17.0 s& 1.10e-04\\
\texttt{FMM3D}\hphantom{abc} \texttt{cpu}   & 8 M & 40    & 33.9 s & 8.19e-05 \\
\texttt{skelFMM cuda}&                            &   320 & 841 ms & 8.30e-05\\
\end{tabular}
\caption{Laplace kernel for points on sphere surface.}
\label{t:comp_lp_sphere}
\end{subfigure}
\caption{\small The figure compares the performance of \texttt{skelFMM}, run on the CPU (using double precision) and GPU (using single precision), to \texttt{FMM3D}, a well-established FMM package implemented in Fortran (using double precision). Only the application time, $T_{\rm apply}$, is reported, noting that \texttt{FMM3D} includes a more efficient precomputation stage. 
}
\label{f:3dlp_comp}
\end{figure}

\vspace{1em}

\input{table_breakdown}

\vspace*{\fill} 

\clearpage

%% file: table_breakdown.tex
\begin{table}[htb!]
\begin{tabular}{|l|rr|rr|rr|rr|}\hline
 & \multicolumn{4}{c|}{\textbf{Sphere} } & \multicolumn{4}{c|}{\textbf{Cube}} \\ 
  & \multicolumn{4}{c|}{Laplace, $N$= 4 M, tol = $10^{-4}$} & \multicolumn{4}{c|}{Laplace, $N$= 4 M, tol = $10^{-5}$} \\ \hline

\textbf{Operation} & \multicolumn{2}{c|}{\texttt{skelFMM cuda}} & \multicolumn{2}{c|}{\texttt{skelFMM cpu}} & \multicolumn{2}{c|}{\texttt{skelFMM cuda}} & \multicolumn{2}{c|}{\texttt{skelFMM cpu}} \\ \hline
upward pass & 60 ms & 10.3 \% & 93 ms & 0.7 \% & 125 ms & 18.7 \% & 120 ms & 0.6 \% \\
add to $\mtx u_{\B}$ (leaf) & 132 ms & 22.8 \% & 6.0 s & 44.4 \% & 106 ms & 15.8 \% & 8.9 s & 43.6 \% \\
add to $\mtx u_{\B}$ (tree) & 104 ms & 17.9 \% & 4.7 s & 34.8 \% & 110 ms & 16.4 \% & 6.5 s & 31.9 \% \\
subtract (ifo) & 89 ms & 15.3 \% & 1.3 s & 9.6 \% & 105 ms & 15.7 \% & 4.6 s & 22.5 \% \\
subtract (tfi,ifs) & 29 ms & 5.0 \% & 800 ms & 5.9 \% & 0 ms & 0.0 \% & 0 ms & 0.0 \% \\
downward pass & 115 ms & 19.8 \% & 150 ms & 1.1 \% & 213 ms & 31.8 \% & 326 ms & 1.6 \% \\\hline
\textbf{Total time} & 580 ms &  & 13.5 s &  & 670 ms & & 20.4 s &  \\\hline
\end{tabular}
\vspace{0.5em}
\caption{\small Detailed breakdown of execution times (in milliseconds or seconds) and their percentage contributions for major operations in the \texttt{skelFMM} algorithm. The analysis highlights the performance benefits of GPU acceleration
 and identifies dominant operations, such as adding to $\mtx{u}_{\B}$ for tree and leaf boxes, as well as suggests
possible room for improvement for more optimized implementations.}
\label{t:3dlp_comp}
\end{table}